\documentclass[11pt]{article}
\usepackage{amsmath,epsfig}
\newcommand{\eps}{\varepsilon}

\newcommand{\Tmrca}{T^{\mathrm{mrca}}}
\newcommand{\Tor}{T^{\mathrm{or}}}
\newcommand{\Nmrca}{N^{\mathrm{mrca}}}
\newcommand{\Nanc}{N^{\mathrm{anc}}}
\newcommand{\sfrac}[2]{{\textstyle\frac{#1}{#2}}}
\newcommand{\ed}{\ \stackrel{d}{=} \ }
\newcommand{\cd}{\ \stackrel{d}{\to} \ }
\newcommand{\cp}{\ \stackrel{p}{\to} \ }

\newcommand{\var}{\mathrm{var}}
\usepackage{amsthm}
\theoremstyle{plain}
\newtheorem{thm}{Theorem}
\newtheorem{lem}[thm]{Lemma}
\newtheorem{cor}[thm]{Corollary}
\newtheorem{prop}[thm]{Proposition}
\theoremstyle{definition}

\theoremstyle{remark}

\newcommand{\T}{\ensuremath{\mathcal{T}}}
\newcommand{\Ttn}{\ensuremath{\mathcal{T}_{t,n}}}
\newcommand{\Ttnn}{\ensuremath{\mathcal{T}_{t_n,n}}}
\newcommand{\TtnnN}{\ensuremath{\mathcal{T}_{t_n,n,N_n}}}
\newcommand{\Tn}{\ensuremath{\mathcal{T}_{n}}}
\newcommand{\Nhat}{\widehat{N}}

\newcommand{\prob}{\mathrm{P}}
\newcommand{\Ex}{\mathrm{E}}
\newcommand{\mean}{\Ex}
\newcommand{\C}{\ensuremath{\mathcal{C}}}

\newcommand{\spt}{\ensuremath{{t}}}
\newcommand{\ptc}{\ensuremath{\mathit{\pi}_{1,\spt}}}
\newcommand{\pt}{\ensuremath{\mathit{\pi}_{1}}}

\newcommand{\A}{\ensuremath{\mathcal{A}}}
\newcommand{\Atn}{\ensuremath{\mathcal{A}_{t,n}}}
\newcommand{\Atnn}{\ensuremath{\mathcal{A}_{t_n,n}}}
\newcommand{\An}{\ensuremath{\mathcal{A}_{n}}}

\begin{document}
\title{{\bf A Critical Branching Process Model for
    Biodiversity}
} 

\author{David Aldous\thanks{Research supported by N.S.F. Grant DMS-0203062}\\ \\
Department of Statistics\\
University of California\\
367 Evans Hall \# 3860\\
Berkeley, CA 94720-3860
\and Lea Popovic\\ \\
IMA\\
University of Minnesota\\
400 Lind Hall\\
Minneapolis, MN 55455}
\maketitle

\begin{abstract}
Motivated as a null model for comparison with data, we
study the following model for a phylogenetic tree on $n$ 
extant species.
The origin of the clade is a random time in the past, whose
(improper) distribution is uniform on $(0,\infty)$. 
After that origin, the process of extinctions and speciations is
a continuous-time critical branching process of constant rate, 
conditioned on having the prescribed number $n$ of species
at the present time.
We study various mathematical properties of this model as
$n \to \infty$ limits:
time of origin and of most recent common ancestor; 
pattern of divergence times within lineage trees; 
time series of numbers of species;
number of extinct species
in total, or ancestral to extant species; and ``local" structure
of the tree itself.
We emphasize several mathematical techniques: associating walks with trees,
a point process representation of lineage trees, and Brownian limits.
\end{abstract}

\newpage
\section{Introduction}
\subsection{A big picture}
This paper forms part of a larger project, which we first outline.
There is a substantial literature on comparing data on different
aspects of {\em biodiversity} or {\em macroevolution} -- the evolutionary history of
speciations and extinctions -- with the predictions of simple
``pure chance" stochastic models.
Available data includes\\
$\bullet$ fossil time series -- fluctuations in number of taxa over time;\\
$\bullet$ shapes of phylogenetic trees on extant species
(Mooers and Heard \cite{mooers} provide an extensive survey);\\
$\bullet$ the distribution of number of species per genus. \\
The fit of simple models, and of more elaborate models 
incorporating conjectured biological process, 
have been studied in these contexts.
While data-motivated models are scientifically natural, 
a mathematical aesthetic suggests a somewhat different approach:
start with a ``pure chance" model which encompasses simultaneously
all the kinds of data that one might hope to find.
Here are two instances of what one would like such a model to provide.\\
$\bullet$ 
Joint description of the phylogenetic tree on an extant clade of species,
its extension to the tree on an observed small proportion of
extinct species, and the (unobserved) entire tree on all
extinct species.\\
$\bullet$ Joint description of fossil time series at different levels of the taxonomic
hierarchy.\\
We emphasize the latter because paleontology
literature tends to assume that a model can be applied at any level,
without enquiring whether this assumption is logically self-consistent.

Our purpose in the larger project is to present what is arguably the mathematically
fundamental such model.
The underlying model at the species level is simple -- a critical branching process conditioned to
have $n$ lineages at the present time.
This is the subject of the present paper,
which is the part of the project most closely related to
classical and contemporary applied probability.
In a subsequent paper aimed at a more biological audience we will describe how
the model extends to higher-order taxa by assuming each new species
has some probability of founding a new higher-order taxon; we will consider
several explicit classification schemes emphasizing desiderata
such as monophyletic groups.

Conceptually, this is a {\em neutral} model which does not 
incorporate conjectured biological process such as
intrinsic tendency for species numbers to increase, differential
speciation or extinction rates, or ecological constraints on numbers
of species.
For well-understood mathematical reasons
(see section \ref{sec:final})
neutral models like ours are implausible for large clades.
In a sense,
the model seems most appropriate as a ``null hypothesis" for
small clades,
at the recent fringe of the Tree of Life,
or for a geological period free of mass extinctions
and their aftermath.

Biological questions motivating our model, and suggested by 
the results of analysis of our model, will be treated in detail
elsewhere, so we give just a brief mention here.
Phylogenetic trees on extant species are nowadays based on
molecular data; technical aspects of tree reconstruction 
form a large and important subject 
\cite{felsenstein-IP,Semple}.
But that is not our focus.  
Let us assume that in the near future we will have a large database
of essentially correct phylogenetic trees, and also assume
these include the time points of divergence of lineages
(rather than just the ``shape" of the tree).
How might one use such a database?
\\
(i) {\em Inference about a particular clade.}
If we have no direct knowledge about extinct species, then we cannot observe 
past fluctuations of number of species with time, and cannot observe
the time of origin of the clade (typically longer ago
than the observable time of most recent common ancestor of extant species).
Inference about such quantities requires some stochastic model;
given a model, one can use the observed phylogenetic tree to make inferences.
\\
(ii) {\em Statistical properties of phylogenetic trees in general.}
In what systematic way do real phylogenetic trees differ from predictions of a simple
model like ours which treats macroevolution as a purely random
process, and what is the biological significance of such differences?

\subsection{Standard models}
Ours is, roughly speaking, the {\em third} simplest model one might devise,
so let us first recall the two simpler models.
\paragraph{The Yule model.}
Yule \cite{yule24} proposed the basic model for speciations
without extinctions.
Initially there is one species.  
Thereafter, independently for each existing species, new species originate
as ``daughter" species at constant rate $\lambda$ 
(i.e. at the times of a Poisson (rate $\lambda$) process).
So for given $n$ one can get a model for an $n$-species tree
by taking the present as a random time at which the number of
species equals $n$.
(The associated continuous-time Markov chain counting number of
species is often called the {\em Yule process}, though its
origin as a model for species is often forgotten.)

\paragraph{The Moran/coalescent model.}
These models, developed and extensively used in population genetics,
can also be applied to macroevolution
(see e.g. \cite{hey92}).
In the Moran model (\cite{ew79} sec. 3.3) the number of 
coexisting species is fixed at $n$.
At successive discrete times, one randomly-chosen species goes extinct
and another randomly chosen one speciates.
Implicit in this model (run from the indefinite past until the present) 
is a model for the phylogenetic tree on the
$n$ extant species;
for large $n$, with suitable rescaling of the time unit, 
the phylogenetic tree approximates the much-studied continuous-time {\em coalescent} model.
To describe the coalescent model, we run time backwards 
from the present, starting with $n$
``lines of descent";
in a time interval $dt$, each {\em pair} of lines of descent has chance $dt$
to merge (``coalesce") into one line, and we continue until reaching
a single {\em most recent common ancestor}.
See \cite{mohle00} for a recent survey.

\paragraph{Why a third model?}
Obviously many basic inference questions mentioned earlier -- about
fluctuations in past numbers of (extinct) species, for instance -- are
not satisfactorily handled within the models above.
Biologists have studied more elaborate models, 
mostly in one of two categories.
(We will give a more detailed account elsewhere, but the bottom
line is that the actual fit of real-world data to parametric models
has not been studied as definitively as one might have expected.)
{\em Exponential growth} models are exemplified by the linear birth-and-death chain
model 
for species numbers
($\lambda_i = \lambda i, \mu_i = \mu i$).
This leads to a model \cite{NMH94} with $3$ parameters 
$(\lambda,\mu,t_*)$
where $t_*$ is time of origin of clade.
{\em Logistic} stochastic models posit a logistic-shaped curve for species numbers,
and also require $3$ or $4$ parameters to specify.
In contrast, 
the model we study (described carefully in Section \ref{sec:model}) has only $1$ parameter
(mean species lifetime).
It is this simplicity, and the desire to avoid the particular biological presumptions
underlying exponential growth or logistic type models, that motivates
our particular model.

\subsection{Outline of results}
A {\em clade} is the set of all species which are descendants of
some (typically extinct) species.
The succinct description 
(to be elaborated in Section \ref{sec:model})
of our model for a phylogenetic tree 
$\Tn$ on a clade with $n$ extant species is as follows.
\begin{quote}
The origin of the clade is a random time in the past, whose
(improper) distribution is uniform on $(0,\infty)$. 
After that origin, the process of extinctions and speciations is
a continuous-time critical branching process of constant rate, 
conditioned on having the prescribed number $n$ of species
at the present time.
\end{quote}
The conditioning is conceptually important: given a real
phylogenetic tree on $23$ extant species, we want to be able to compare
it to the predictions of a stochastic model generating trees
on exactly $23$ extant species.
The uniform prior (for time of origin) avoids the necessity to introduce a second
parameter into the model.

There is a vast mathematical literature on branching processes,
but we haven't found detailed discussion of any very similar model.
In the biological literature, Wollenberg et al \cite{WAA96} give
a simulation study of a similar model.
On the other hand, this model is clearly open to analysis by
the known techniques of applied probability.
We exploit one particular modern approach to classical branching process theory:
representing trees as walks.
See Pitman \cite{csp} for a recent survey.
This both leads to an exact ``point process" description of lineage tree distributions
(Proposition \ref{prop:cond_ptpr})
and permits us to study asymptotics via weak convergence to Brownian motion.
This methodology is known to specialists in other aspects
of random trees but is perhaps less familiar in the subject
of phylogenetic trees, so we try to explain the key ideas carefully
even though they are not entirely new.

Our results describe distributional properties of
various aspects of the tree $\T_n$.
\begin{itemize}
\item The lineage tree, via exact formulas (Proposition \ref{prop:cond_ptpr}),
``global limits" (Corollaries \ref{cor:conv_ptpr} and \ref{cor:global}), and
``local limits" (Corollaries \ref{cor:local} and \ref{calc}).
\item The time series of number of species (Lemma \ref{lem:timerev}),
the maximum number of coexisting species (Corollary \ref{cor:max_popn}),
and the total number of extinct species (Corollary \ref{cor:conv_popn}).
\item The local limit structure of the complete (i.e. including
extinct species) phylogenetic tree, relative to either a typical
extant species (Proposition \ref{thm5}) or a typical extinct species (Proposition \ref{thm4}).
\item The joint distribution of time of origin of clade and
time of most recent common ancestor(Corollary \ref{cor:t_or,mrca}), joint
also with the number of species alive at the time of most recent
common ancestor (Corollary \ref{cor:conv_triple}).
\item The number of extinct species ancestral to some extant species 
(Corollary \ref{cor:conv_anc}).
\end{itemize}

Finally we should admit that the whole paradigm of studying
$n \to \infty$ asymptotics is rather unnatural, because the model
is biologically unrealistic for large $n$, but 
one can hope that the approximations implicit in asymptotic results
are qualitatively correct for smaller values of $n$.
Our web site \cite{me-Phylo-site} shows Monte Carlo 
simulations for $n = 8, 12, 20$ with $10$ repetitions.
One can check that numerical values are broadly consistent with 
the asymptotic predictions.

\section{Model and notation}\label{sec:model}

In stating and deriving mathematical results 
we use the traditional language of branching processes 
({\em individuals, children, births, \ldots})
even though we are envisaging species 
and so should be writing 
({\em species, daughter species, speciations, \ldots}).

Let $\T$ be a continuous time critical branching process (CBP)
starting with one individual. According to this process each
individual lives for an Exponential (rate $\lambda$) time, for some
$\lambda>0$,  during
which it gives birth at times of an independent Poisson (rate
$\lambda$) process. After birth all individuals behave independently
of each other. 
We can and will scale time so that $\lambda = 1$; so the time unit
is interpreted as mean species lifetime.

Write
$N_\T(t) \geq 0$ for the number of individuals alive at time
$t$ after the origin of $\T$.
A classical result
(\cite{fel68} \S XVII.10.11) 
gives a modified Geometric distribution
\begin{eqnarray}
\prob(N_\T(t) = n ) &=& \frac{t^{n-1}}{(1+t)^{n+1}}, \ n \geq 1 \label{Feller}\\
&=& \frac{t}{1+t}, \ n = 0 . \nonumber
\end{eqnarray}
Write $\Ttn$ for the process $\T$ originating 
at time $t$ in the past and conditioned on having exactly $n$
individuals at the present time.
Within a process like $\Ttn$ or $\Tn$ below, we use the
notational convention that 
``time $s$" means time $s$ before present.
Thus within $\Ttn$, the time parameter $s$ decreases from $t$ to $0$,
meaning that time increases from time $t$ before present to 
time $0$, the present time.
Our previous verbal definition of our model $\Tn$ 
as a Bayes posterior (for $\T$ started at a uniform past time 
$t \in (0,\infty)$ and conditioned on having $n$ individuals
at the present time $0$) 
now becomes the following rigorous definition.
Fix $n \geq 1$.
\[
\prob(\Tn \in \cdot) = 
\frac{\int_0^\infty 
\prob(\Ttn \in \cdot) 
\prob(N_\T(t) = n) 
\ dt }
{\int_0^\infty 
\prob(N_\T(t) = n) 
\ dt }
. \]
Using (\ref{Feller}) and the calculus fact
\[
\int_0^\infty\frac{s^{n-1}}{(1+s)^{n+1}}\ ds
= n^{-1} 
\]
turns this into
\begin{equation}
\prob(\Tn \in \cdot) = 
\int_0^\infty 
\prob(\Ttn \in \cdot) 
\ \sfrac{nt^{n-1}}{(1+t)^{n+1}} \ dt . \label{Tmix}
\end{equation}
Within the random tree $\Tn$, the time parameter $s$ 
decreases from a random ``time of origin" 
$\Tor_n$ to $0$, where by 
the formula above,
$\Tor_n$ has density function
\begin{equation}
q_n(t)=\frac{nt^{n-1}}{(1+t)^{n+1}}, \quad t>0
. \label{def:qnt}
\end{equation}
We shall refer to $\Tn$ and $\Ttn$ as the {\em complete trees}. 
Returning to biological terminology, a complete tree 
records the birth times and every (extinct or extant) species in
a clade and the extinction time of extinct species.
Every
realization of 
a complete tree  also uniquely determines a realization of a {\em
  lineage tree}  of the extant species. This is the
smallest subtree of the complete tree that contains all the divergence
times for pairs of lineages 
of extant species, without recording which ancestral species contain the lineage. We let $\Atn$ and $\An$ denote the lineage trees
of $\Ttn$ and $\Tn$ respectively.  
The time parameter $s$ within $\An$ decreases from the
time $\Tmrca_n$ of 
{\em most recent common ancestor} 
of the $n$ extant species, to the present time $0$.
(The lineage tree is what is usually called the 
{\em phylogenetic tree}, though logically all the trees
under consideration are different kinds of phylogenetic tree.)

\section{Point process representations of lineage trees}\label{sec:loc_anc}
\subsection{An exact description}
It is perhaps remarkable that there is a useful exact description
of the lineage tree $\Atn$, based on a certain
{\em point process representation}
illustrated in Figure 1.
Consider an arbitrary lineage tree on $n$ species.
Draw the tree as in Figure 1; 
recursively from the top down, at each divergence point of lineages 
choose randomly which branch is drawn on the left and which on the right.
After drawing the tree, label species as 
$1,2,\ldots,n$ in left-to-right order.
Each divergence of lineages involves adjacent contiguous blocks of species,
say 
$\{i,i+1,\ldots,j\}$ and $\{j+1,j+2,\ldots,k\}$,
and occurs at some time $s$.
We mark the occurrence of this divergence by a mark
$\times$ at coordinates $(j+ \frac{1}{2},s)$
and then draw the combined lineage as a vertical line upwards from the mark.

\setlength{\unitlength}{0.25in}
\begin{picture}(22,11)(2,-1)
\put(11,0){\line(0,1){9}}
\put(10.4,-0.17){0}
\put(10.4,1.83){2}
\put(10.4,3.83){4}
\put(9.8,5){time}
\put(10.4,5.83){6}
\put(10.4,7.83){8}
\put(11,0){\line(-1,0){0.2}}
\put(11,2){\line(-1,0){0.2}}
\put(11,4){\line(-1,0){0.2}}
\put(11,6){\line(-1,0){0.2}}
\put(11,8){\line(-1,0){0.2}}
\put(1.27,0.24){$\times$}
\put(2.27,2.31){$\times$}
\put(3.27,0.59){$\times$}
\put(4.27,4.75){$\times$}
\put(5.27,0.11){$\times$}
\put(6.27,0.49){$\times$}
\put(7.27,1.59){$\times$}
\put(8.27,7.52){$\times$}
\put(9.27,1.24){$\times$}
\put(1,0){\line(0,1){0.37}}
\put(2,0){\line(0,1){0.37}}
\put(1,0.37){\line(1,0){1}}
\put(1.5,2.44){\line(0,-1){2.07}}
\put(1.5,2.44){\line(1,0){2}}
\put(3,0){\line(0,1){0.72}}
\put(4,0){\line(0,1){0.72}}
\put(3,0.72){\line(1,0){1}}
\put(3.5,0.72){\line(0,1){1.72}}
\put(5,0){\line(0,1){0.24}}
\put(6,0){\line(0,1){0.24}}
\put(5,0.24){\line(1,0){1}}
\put(7,0){\line(0,1){0.62}}
\put(5.5,0.62){\line(0,-1){0.38}}
\put(5.5,0.62){\line(1,0){1.5}}
\put(8,0){\line(0,1){1.72}}
\put(6.5,1.72){\line(0,-1){1.1}}
\put(6.5,1.71){\line(1,0){1.5}}
\put(2.5,4.88){\line(0,-1){2.44}}
\put(2.5,4.88){\line(1,0){5}}
\put(7.5,1.72){\line(0,1){3.16}}
\put(9,0){\line(0,1){1.37}}
\put(10,0){\line(0,1){1.37}}
\put(9,1.37){\line(1,0){1}}
\put(9.5,1.37){\line(0,1){6.28}}
\put(4.5,4.88){\line(0,1){2.77}}
\put(4.5,7.65){\line(1,0){5}}
\put(8.5,7.65){\line(0,1){1.35}}
\put(12.27,0.24){$\times$}
\put(13.27,2.31){$\times$}
\put(14.27,0.59){$\times$}
\put(15.27,4.75){$\times$}
\put(16.27,0.11){$\times$}
\put(17.27,0.49){$\times$}
\put(18.27,1.59){$\times$}
\put(19.27,7.52){$\times$}
\put(20.27,1.24){$\times$}
\put(11,0){\line(1,0){10}}
\put(21,0){\line(0,1){9}}
\put(11,0){\line(0,-1){0.2}}
\put(13,0){\line(0,-1){0.2}}
\put(15,0){\line(0,-1){0.2}}
\put(17,0){\line(0,-1){0.2}}
\put(19,0){\line(0,-1){0.2}}
\put(21,0){\line(0,-1){0.2}}
\put(10.83,-0.8){0}
\put(12.83,-0.8){2}
\put(14.83,-0.8){4}
\put(16.83,-0.8){6}
\put(18.83,-0.8){8}
\put(20.73,-0.8){10}
\end{picture}

{\bf Figure 1.}
{\small 
The point process representation of a lineage tree on $n = 10$ species.}

\vspace{0.12in}
\noindent
The advantage of this precise way of drawing the tree is that one can clearly
reconstruct the tree from the coordinates 
$\{(i+\frac{1}{2},s_i), 1 \leq i \leq n-1\}$
of the marks.
So the distribution of the point process of marks 
serves to specify the distribution of the lineage tree.
\begin{prop}[\cite{lea-paper} Lemma 3]
\label{prop:cond_ptpr}
Fix $n \geq 2$ and $t > 0$.
The point process
$\{(i+\frac{1}{2},h_i), \ 1 \leq i \leq n-1\}$
where the $(h_i)$ are i.i.d. with density function
\begin{equation} 
f_t(s) = (1+t^{-1})(1+s)^{-2}, \ 0<s<t 
\label{def:fts}
\end{equation}
represents the lineage tree $\Atn$ within the complete tree $\Ttn$.
\end{prop}
The derivation of this result will be explained in Section \ref{sec:contour},
where the underlying {\em contour process} is exploited further.

We are mostly concerned with the lineage tree $\An$,
which by (\ref{Tmix}) has a mixture representation
\begin{equation}
\prob(\An \in \cdot) = \int_0^\infty \prob(\Atn \in \cdot) q_n(t) \ dt 
\label{Amix}
\end{equation}
where $q_n(t)$ is the density function (\ref{def:qnt}) of $\Tor_n$.
One can get exact formulas for various attributes of $\An$.
Consider for instance the number of lineages at time $s$.
Because each divergence time creates one extra lineage, it is clear
that within $\Atn$ this number of lineages is distributed as
\[ 1 + {\rm Binomial}(n-1,\bar{F}_t(s)) \]
for
\[
\bar{F}_t(s) = \int_s^t f_t(u) \ du = 
\frac{t-s}{t(1+s)} . \]
Thus within $\An$ the distribution is the mixture of Binomials implied 
by (\ref{Amix}).
Similarly, the exact distribution of the time 
$\Tmrca_n$ of most recent common ancestor
is
\[ \prob(\Tmrca_n \leq u) = 
\int_0^\infty \left(1 - \bar{F}_t(u) \right)^{n-1} \ q_n(t) \ dt . \]
In this paper we focus on $n \to \infty$ asymptotics, which may give more
conceptual insight than do complicated exact formulas.
As we see below, it is useful to distinguish two kinds of asymptotics:
{\em global limits} refer to times of order $n$, whereas
{\em local limits} refer to times of order $1$.

\subsection{The global limit point process}
From the formula (\ref{def:qnt}) for $q_n(t)$ we calculate: 
if $t_n/n \to t>0$ then
\[nq_n(t_n)=\frac{n^2}{(1+t_n)^2}\left(1-\frac{1}{1+t_n}\right)^{n-1}
  \mathop{\longrightarrow}_{n\to\infty}t^{-2}\,e^{-\frac{1}{t}}
.\]
The limit is the density function of the
{\em Inverse Exponential} 
IE(1) distribution, that is to say of
$1/\xi$ where $\xi$ has Exponential(1) distribution.
So we have shown
\begin{lem}
\label{lem:t_or}
$n^{-1} \Tor_n \cd \Tor$, say, where the limit
$\Tor$ has IE(1) distribution.
\end{lem}
Now reconsider Figure 1.
To obtain a global limit we want to rescale time by a factor $n$
and we want to rescale the left-to-right positions of marks to fit into
a unit interval $[0,1]$, implying they also must be rescaled by a factor $n$.
Thus the original point process of marks
$\{(i+\frac{1}{2},s_i), 1 \leq i \leq n-1\}$
is rescaled to
$\{(\sfrac{i+\frac{1}{2}}{n},\sfrac{s_i}{n}), 1 \leq i \leq n-1\}$.
In the setting of Proposition \ref{prop:cond_ptpr},
the relevant calculation is:
\[ \mbox{
if $s_n/n \to s> 0$ and $t_n/n \to t>0$ then
$n^2f_{t_n}(s_n) \to s^{-2}$} \]
and the following limit behavior is intuitively clear.
\begin{cor}[\cite{lea-paper} Lemma 4, Theorem 5]
\label{cor:conv_ptpr}
Let $t_n/n \to t > 0$.
The rescaled point process
$\{(\sfrac{i+\frac{1}{2}}{n},\sfrac{h_i}{n}), 1 \leq i \leq n-1\}$
associated with the lineage tree $\Atnn$ converges in distribution
to the Poisson point process 
($\ptc$, say) whose intensity measure is 
  $\nu(dl\times ds)=dl\,s^{-2}ds  \; \mathbf{1}_{[0,1]\times (0,t)}.$
\end{cor}
The limit $\ptc$, illustrated in Figure 2, has an infinite number of
points close to the lower boundary, but weak convergence on the
open interval $(0,t)$ means convergence over regions away from this boundary.
Figure 2 indicates visually how the Poisson point process limit defines a limit random tree 
which is a kind of ``continuum tree" with a lineage for each real $l \in (0,1)$,
though we do not seek to formalize this idea.

\setlength{\unitlength}{.25in}
\begin{picture}(22,12)(2,-1)
\put(11,0){\line(0,1){10}}
\put(10.4,-0.12){0}
\put(10.1,1.83){10}
\put(10.1,3.83){20}
\put(9.8,5){time}
\put(10.1,5.83){30}
\put(10.1,7.83){40}
\put(10.1,9.83){50}
\put(11,0.15){\line(-1,0){0.2}}
\put(11,2){\line(-1,0){0.2}}
\put(11,4){\line(-1,0){0.2}}
\put(11,6){\line(-1,0){0.2}}
\put(11,8){\line(-1,0){0.2}}

\put(0.08,0.3249){\tiny{$\times$}}
\put(0.10,0.4751){\tiny{$\times$}}
\put(0.38,0.3027){\tiny{$\times$}}
\put(0.51,0.6494){\tiny{$\times$}}
\put(0.52,0.2473){\tiny{$\times$}}
\put(0.53,0.3556){\tiny{$\times$}}
\put(0.55,0.3111){\tiny{$\times$}}
\put(0.65,0.8721){\tiny{$\times$}}
\put(0.72,0.4549){\tiny{$\times$}}
\put(0.67,0.4497){\tiny{$\times$}}

\put(0.69,0.2324){\tiny{$\times$}}
\put(0.72,0.2648){\tiny{$\times$}}
\put(0.73,0.2683){\tiny{$\times$}}
\put(0.88,0.3403){\tiny{$\times$}}
\put(0.89,0.2899){\tiny{$\times$}}
\put(1.11,0.2007){\tiny{$\times$}}
\put(1.21,2.6858){\tiny{$\times$}}
\put(1.30,0.4026){\tiny{$\times$}}
\put(1.41,0.2481){\tiny{$\times$}}
\put(1.43,0.2113){\tiny{$\times$}}

\put(1.92,0.6989){\tiny{$\times$}}
\put(1.96,1.7903){\tiny{$\times$}}
\put(2.06,0.2996){\tiny{$\times$}}
\put(2.09,0.6615){\tiny{$\times$}}
\put(2.12,0.3097){\tiny{$\times$}}
\put(2.12,0.5036){\tiny{$\times$}}
\put(2.26,0.2279){\tiny{$\times$}}
\put(2.56,9.4492){\tiny{$\times$}}
\put(2.60,0.6986){\tiny{$\times$}}
\put(2.63,0.5583){\tiny{$\times$}}

\put(2.74,0.3402){\tiny{$\times$}}
\put(2.75,0.2974){\tiny{$\times$}}
\put(2.77,0.5594){\tiny{$\times$}}
\put(2.83,0.5048){\tiny{$\times$}}
\put(2.84,1.0440){\tiny{$\times$}}
\put(2.86,0.6921){\tiny{$\times$}}
\put(3.03,0.3171){\tiny{$\times$}}
\put(3.23,0.8525){\tiny{$\times$}}
\put(3.30,0.2156){\tiny{$\times$}}
\put(3.37,0.5568){\tiny{$\times$}}

\put(3.66,0.4625){\tiny{$\times$}}
\put(3.68,0.4199){\tiny{$\times$}}
\put(3.71,0.4387){\tiny{$\times$}}
\put(3.76,0.3442){\tiny{$\times$}}
\put(3.77,0.4563){\tiny{$\times$}}
\put(4.37,0.8363){\tiny{$\times$}}
\put(4.50,0.2456){\tiny{$\times$}}
\put(4.61,0.4595){\tiny{$\times$}}
\put(4.95,0.8953){\tiny{$\times$}}
\put(5.00,1.8198){\tiny{$\times$}}

\put(5.15,0.5579){\tiny{$\times$}}
\put(5.23,0.3461){\tiny{$\times$}}
\put(5.36,0.2359){\tiny{$\times$}}
\put(5.56,0.2197){\tiny{$\times$}}
\put(5.60,0.2449){\tiny{$\times$}}
\put(5.72,1.1909){\tiny{$\times$}}
\put(5.74,0.2539){\tiny{$\times$}}
\put(5.77,0.3690){\tiny{$\times$}}
\put(5.80,0.6208){\tiny{$\times$}}
\put(5.81,0.2964){\tiny{$\times$}}

\put(5.91,0.3886){\tiny{$\times$}}
\put(5.92,0.2789){\tiny{$\times$}}
\put(6.39,2.2199){\tiny{$\times$}}
\put(6.40,0.8050){\tiny{$\times$}}
\put(6.46,0.2718){\tiny{$\times$}}
\put(6.50,0.3196){\tiny{$\times$}}
\put(6.61,0.2672){\tiny{$\times$}}
\put(6.62,1.4153){\tiny{$\times$}}
\put(6.71,1.4675){\tiny{$\times$}}
\put(6.79,1.9069){\tiny{$\times$}}

\put(6.93,0.5382){\tiny{$\times$}}
\put(6.97,0.7191){\tiny{$\times$}}
\put(7.37,4.2978){\tiny{$\times$}}
\put(7.44,0.5524){\tiny{$\times$}}
\put(7.62,1.2070){\tiny{$\times$}}
\put(7.63,0.2541){\tiny{$\times$}}
\put(7.68,2.5074){\tiny{$\times$}}
\put(8.07,0.2109){\tiny{$\times$}}
\put(8.32,0.2692){\tiny{$\times$}}
\put(8.32,1.9950){\tiny{$\times$}}

\put(8.37,0.3845){\tiny{$\times$}}
\put(8.40,0.2038){\tiny{$\times$}}
\put(8.47,0.2075){\tiny{$\times$}}
\put(8.49,3.8511){\tiny{$\times$}}
\put(8.54,0.2262){\tiny{$\times$}}
\put(8.55,0.6454){\tiny{$\times$}}
\put(8.63,0.2315){\tiny{$\times$}}
\put(8.67,0.4701){\tiny{$\times$}}
\put(8.72,1.1079){\tiny{$\times$}}
\put(8.76,0.2102){\tiny{$\times$}}

\put(9.00,0.3184){\tiny{$\times$}}
\put(9.05,0.6098){\tiny{$\times$}}
\put(9.23,0.2230){\tiny{$\times$}}
\put(9.29,0.3926){\tiny{$\times$}}
\put(9.33,0.8723){\tiny{$\times$}}
\put(9.38,0.2772){\tiny{$\times$}}
\put(9.40,0.2674){\tiny{$\times$}}
\put(9.48,0.3766){\tiny{$\times$}}
\put(9.61,0.2408){\tiny{$\times$}}
\put(9.64,0.5911){\tiny{$\times$}}

\put(2.71,9.5492){\line(0,1){0.4}}
\put(1.36,2.7858){\line(0,1){6.7634}}
\put(7.52,4.3978){\line(0,1){5.1514}}
\put(0.8,0.9721){\line(0,1){1.8137}}
\put(2.11,1.8903){\line(0,1){0.8955}}
\put(6.54,2.3199){\line(0,1){2.0779}}
\put(8.64,3.9511){\line(0,1){0.4467}}
\put(0.66,0.7494){\line(0,1){0.2227}}
\put(0.87,0.5549){\line(0,1){0.4172}}
\put(2.07,0.7989){\line(0,1){1.0914}}
\put(2.24,0.7615){\line(0,1){1.1288}}
\put(5.15,1.9198){\line(0,1){0.4001}}
\put(6.94,2.0069){\line(0,1){0.313}}
\put(7.83,2.6074){\line(0,1){1.3437}}
\put(8.87,1.2079){\line(0,1){2.7432}}
\put(2.99,1.144){\line(0,1){0.7758}}
\put(5.87,1.2909){\line(0,1){0.628}}
\put(6.86,1.5675){\line(0,1){0.4394}}
\put(7.12,0.8191){\line(0,1){1.1878}}
\put(7.77,1.307){\line(0,1){1.3004}}
\put(8.49,2.095){\line(0,1){0.5124}}
\put(8.7,0.7454){\line(0,1){0.4625}}
\put(9.48,0.9723){\line(0,1){0.2356}}

\put(1.36,9.5492){\line(1,0){6.16}}
\put(0.8,2.7858){\line(1,0){1.31}}
\put(6.54,4.3978){\line(1,0){2.0812}}
\put(0.66,0.9721){\line(1,0){0.22}}
\put(2.07,1.8903){\line(1,0){0.17}}
\put(5.15,2.3199){\line(1,0){1.79}}
\put(7.83,3.9511){\line(1,0){1.04}}
\put(0.25,0.7494){\line(1,0){0.55}}
\put(2.99,1.9198){\line(1,0){2.88}}
\put(8.7,1.2079){\line(1,0){0.78}}
\put(7.77,2.6074){\line(1,0){0.72}}
\put(6.86,2.0069){\line(1,0){0.26}}
\put(5.1,1.2909){\line(1,0){1.45}}
\put(0.25,0.5751){\line(0,1){0.1743}}
\put(0.8,0.7494){\line(0,-1){0.25}}
\put(5.1,0.9953){\line(0,1){0.2956}}
\put(6.55,0.905){\line(0,1){0.3859}}
\put(8.24,2.095){\line(0,-1){1.8}}
\put(8.56,2.095){\line(0,-1){1.6}}
\put(8.24,2.095){\line(1,0){0.31}}
\put(9.18,0.9723){\line(1,0){0.6}}
\put(9.18,0.9723){\line(0,-1){0.32}}
\put(9.78,0.9723){\line(0,-1){0.35}}
\put(7.61,1.307){\line(1,0){0.2}}
\put(7.61,1.307){\line(0,-1){0.65}}
\put(7.8,1.307){\line(0,-1){0.95}}
\put(6.73,1.5675){\line(1,0){0.3}}
\put(6.73,1.5675){\line(0,-1){1.25}}
\put(7.04,1.5675){\line(0,-1){0.95}}
\put(5.98,0.905){\line(1,0){0.7}}
\put(5.98,0.905){\line(0,-1){0.25}}
\put(6.68,0.905){\line(0,-1){0.45}}
\put(4.5,0.99){\line(1,0){0.8}}
\put(4.5,0.99){\line(0,-1){0.1}}
\put(5.3,0.99){\line(0,-1){0.35}}
\put(3.92,0.91){\line(1,0){0.81}}
\put(3.92,0.91){\line(0,-1){0.29}}
\put(4.75,0.91){\line(0,-1){0.34}}
\put(2.74,1.144){\line(1,0){0.65}}
\put(2.74,1.144){\line(0,-1){0.3}}
\put(3.39,1.144){\line(0,-1){0.25}}
\put(3.19,0.994){\line(1,0){0.36}}
\put(3.19,0.994){\line(0,-1){0.54}}
\put(3.55,0.994){\line(0,-1){0.29}}
\put(1.44,0.7989){\line(1,0){0.79}}
\put(1.44,0.7989){\line(0,-1){0.29}}
\put(2.21,0.7989){\line(0,-1){0.39}}

\put(11.08,0.3249){\tiny{$\times$}}
\put(11.10,0.4751){\tiny{$\times$}}
\put(11.38,0.3027){\tiny{$\times$}}
\put(11.51,0.6494){\tiny{$\times$}}
\put(11.52,0.2473){\tiny{$\times$}}
\put(11.53,0.3556){\tiny{$\times$}}
\put(11.55,0.3111){\tiny{$\times$}}
\put(11.65,0.8721){\tiny{$\times$}}
\put(11.72,0.4549){\tiny{$\times$}}
\put(11.67,0.4497){\tiny{$\times$}}
\put(11.69,0.2324){\tiny{$\times$}}
\put(11.72,0.2648){\tiny{$\times$}}
\put(11.73,0.2683){\tiny{$\times$}}
\put(11.88,0.3403){\tiny{$\times$}}
\put(11.89,0.2899){\tiny{$\times$}}
\put(12.11,0.2007){\tiny{$\times$}}
\put(12.21,2.6858){\tiny{$\times$}}
\put(12.30,0.4026){\tiny{$\times$}}
\put(12.41,0.2481){\tiny{$\times$}}
\put(12.43,0.2113){\tiny{$\times$}}
\put(12.92,0.6989){\tiny{$\times$}}
\put(12.96,1.7903){\tiny{$\times$}}
\put(13.06,0.2996){\tiny{$\times$}}
\put(13.09,0.6615){\tiny{$\times$}}
\put(13.12,0.3097){\tiny{$\times$}}
\put(13.12,0.5036){\tiny{$\times$}}
\put(13.26,0.2279){\tiny{$\times$}}
\put(13.56,9.4492){\tiny{$\times$}}
\put(13.60,0.6986){\tiny{$\times$}}
\put(13.63,0.5583){\tiny{$\times$}}
\put(13.74,0.3402){\tiny{$\times$}}
\put(13.75,0.2974){\tiny{$\times$}}
\put(13.77,0.5594){\tiny{$\times$}}
\put(13.83,0.5048){\tiny{$\times$}}
\put(13.84,1.0440){\tiny{$\times$}}
\put(13.86,0.6921){\tiny{$\times$}}
\put(14.03,0.3171){\tiny{$\times$}}
\put(14.23,0.8525){\tiny{$\times$}}
\put(14.30,0.2156){\tiny{$\times$}}
\put(14.37,0.5568){\tiny{$\times$}}
\put(14.66,0.4625){\tiny{$\times$}}
\put(14.68,0.4199){\tiny{$\times$}}
\put(14.71,0.4387){\tiny{$\times$}}
\put(14.76,0.3442){\tiny{$\times$}}
\put(14.77,0.4563){\tiny{$\times$}}
\put(15.37,0.8363){\tiny{$\times$}}
\put(15.50,0.2456){\tiny{$\times$}}
\put(15.61,0.4595){\tiny{$\times$}}
\put(15.95,0.8953){\tiny{$\times$}}
\put(16.00,1.8198){\tiny{$\times$}}
\put(16.15,0.5579){\tiny{$\times$}}
\put(16.23,0.3461){\tiny{$\times$}}
\put(16.36,0.2359){\tiny{$\times$}}
\put(16.56,0.2197){\tiny{$\times$}}
\put(16.60,0.2449){\tiny{$\times$}}
\put(16.72,1.1909){\tiny{$\times$}}
\put(16.74,0.2539){\tiny{$\times$}}
\put(16.77,0.3690){\tiny{$\times$}}
\put(16.80,0.6208){\tiny{$\times$}}
\put(16.81,0.2964){\tiny{$\times$}}
\put(16.91,0.3886){\tiny{$\times$}}
\put(16.92,0.2789){\tiny{$\times$}}
\put(17.39,2.2199){\tiny{$\times$}}
\put(17.40,0.8050){\tiny{$\times$}}
\put(17.46,0.2718){\tiny{$\times$}}
\put(17.50,0.3196){\tiny{$\times$}}
\put(17.61,0.2672){\tiny{$\times$}}
\put(17.62,1.4153){\tiny{$\times$}}
\put(17.71,1.4675){\tiny{$\times$}}
\put(17.79,1.9069){\tiny{$\times$}}
\put(17.93,0.5382){\tiny{$\times$}}
\put(17.97,0.7191){\tiny{$\times$}}
\put(18.37,4.2978){\tiny{$\times$}}
\put(18.44,0.5524){\tiny{$\times$}}
\put(18.62,1.2070){\tiny{$\times$}}
\put(18.63,0.2541){\tiny{$\times$}}
\put(18.68,2.5074){\tiny{$\times$}}
\put(19.07,0.2109){\tiny{$\times$}}
\put(19.32,0.2692){\tiny{$\times$}}
\put(19.32,1.9950){\tiny{$\times$}}
\put(19.37,0.3845){\tiny{$\times$}}
\put(19.40,0.2038){\tiny{$\times$}}
\put(19.47,0.2075){\tiny{$\times$}}
\put(19.49,3.8511){\tiny{$\times$}}
\put(19.54,0.2262){\tiny{$\times$}}
\put(19.55,0.6454){\tiny{$\times$}}
\put(19.63,0.2315){\tiny{$\times$}}
\put(19.67,0.4701){\tiny{$\times$}}
\put(19.72,1.1079){\tiny{$\times$}}
\put(19.76,0.2102){\tiny{$\times$}}
\put(20.00,0.3184){\tiny{$\times$}}
\put(20.05,0.6098){\tiny{$\times$}}
\put(20.23,0.2230){\tiny{$\times$}}
\put(20.29,0.3926){\tiny{$\times$}}
\put(20.33,0.8723){\tiny{$\times$}}
\put(20.38,0.2772){\tiny{$\times$}}
\put(20.40,0.2674){\tiny{$\times$}}
\put(20.48,0.3766){\tiny{$\times$}}
\put(20.61,0.2408){\tiny{$\times$}}
\put(20.64,0.5911){\tiny{$\times$}}
\put(11,0.15){\line(1,0){10}}
\put(21,0.15){\line(0,1){10}}
\put(11,0.15){\line(0,-1){0.2}}
\put(13,0.15){\line(0,-1){0.2}}
\put(15,0.15){\line(0,-1){0.2}}
\put(17,0.15){\line(0,-1){0.2}}
\put(19,0.15){\line(0,-1){0.2}}
\put(21,0.15){\line(0,-1){0.2}}
\put(10.83,-0.65){0}
\put(12.83,-0.65){20}
\put(14.83,-0.65){40}
\put(16.83,-0.65){60}
\put(18.83,-0.65){80}
\put(20.73,-0.65){100}
\end{picture}

{\bf Figure 2.}
{\small 
The point process $\pi_{1,t}$ on the right represents the lineage tree
of a continuum of species on the left.}

\vspace{0.12in}
\noindent
The mixture representation (\ref{Amix}) and Corollary \ref{cor:conv_ptpr}
immediately imply a global limit theorem for $\An$.
To state it, let $\Tor$ have IE(1) distribution.
Define a Cox point process 
($\pt$, say)
on $(0,1) \times (0,\infty)$
as follows.
Given $\Tor = t$,
 let $\pt$ be a Poisson point process with the law of $\ptc$.
\begin{cor}
\label{cor:global}
The rescaled point process
$\{(\sfrac{i+\frac{1}{2}}{n},\sfrac{s_i}{n}), 1 \leq i \leq n-1\}$
associated with the lineage tree $\An$, considered jointly with $\Tor_n$, converges in distribution
to the Cox point process 
$\pt$, considered jointly with $\Tor$.
\end{cor}

Here is a quick application of this global limit theorem.
\begin{cor}
\label{cor:t_or,mrca}
The limit joint behavior of 
$\Tor_n$ and $\Tmrca_n$ is given by
\[ (n^{-1}\Tor_n,n^{-1}\Tmrca_n) \cd
(\Tor,\Tmrca) \]
where the limit law has joint density
\[
f_{\Tor,\Tmrca}(t,s)=t^{-2}s^{-2}e^{-\frac{1}{s}},\quad 0<s<t.
\]
The marginal density of $\Tmrca$ is
  \[f_{\Tmrca}(s)=s^{-3}e^{-\frac{1}{s}}, \quad s>0.\]
The limit joint distribution can alternatively be expressed as
  $(\Tor,\Tmrca)\ed
  (\frac{1}{\xi_1},\frac{1}{\xi_1+\xi_2})$
  where $\xi_1,\xi_2$ are i.i.d. Exponential$(1)$.\\
\end{cor}
\begin{proof}
Corollary \ref{cor:global} implies convergence in distribution
to the limit $(\Tor,\Tmrca)$ in which
$\Tmrca$ is defined as the maximum height
(that is, maximum second coordinate)
of any point of $\pt$.
Given that $\Tor=\spt$, $\pt$ is distributed as a Poisson
point process $\ptc$ with intensity measure  
$\nu(dl\times ds)=dl\,s^{-2}ds\; \mathbf{1}_{[0,1]\times (0,t)}$. 
Consequently, for
the conditional law of $\Tmrca$ given $\Tor=t$ we have
\begin{eqnarray*}
\prob(\Tmrca\leq s|\Tor=t)
&=& \prob\big(\{\ptc\cap[0,1]\times(s,\spt)\}=\emptyset\big)
\\
&=& \exp\left(- \int_s^t u^{-2} du \right) \\
&=& e^{\frac{1}{\spt}-\frac{1}{s}} ,\quad 0<s<\spt.
\end{eqnarray*}
So
\[ \prob(\Tmrca \leq s, \Tor \in dt)
= e^{\frac{1}{t}-\frac{1}{s}} P(\Tor \in dt) =
 t^{-2} e^{-\frac{1}{s}} \ dt, \quad 0<s<\spt \]
implying the formula for joint density.
The remaining calculations are straightforward.
\end{proof}

\subsection{The local limit point process}
\label{sec:llp}
There is a different limit regime in which time is not rescaled.
This tells us the local structure of the lineage tree relative
to a given typical species, where ``local" refers to lineages merging
with the given lineage within bounded time.
The relevant calculation is that,
in the setting of Proposition 
\ref{prop:cond_ptpr}, 
if $t_n \to \infty$ then 
\[
f_{t_n}(s) \to f(s):= (1+s)^{-2}, \quad 
0<s<\infty .
\]
Consider the point process on
$(\mathbf{Z} + \frac{1}{2}) \times (0,\infty)$
consisting of points
$\{(i+\frac{1}{2},\eta_i), \ i \in \mathbf{Z}\}$
for i.i.d. $(\eta_i)$ with density
$f(s) = (1+s)^{-2}$.
As illustrated in Figure 3, the point process defines
an infinite tree, $\A_\infty$ say,
on lineages labeled by $\mathbf{Z}$.

\setlength{\unitlength}{0.25in}
\begin{picture}(11,12)
\put(0,0){\line(0,1){11}}
\put(-0.6,-0.17){0}
\put(-0.6,1.83){2}
\put(-0.6,3.83){4}
\put(-1.8,5){time}
\put(-0.6,5.83){6}
\put(-0.6,7.83){8}
\put(-0.8,9.83){10}
\put(0,0){\line(-1,0){0.2}}
\put(0,2){\line(-1,0){0.2}}
\put(0,4){\line(-1,0){0.2}}
\put(0,6){\line(-1,0){0.2}}
\put(0,8){\line(-1,0){0.2}}
\put(0,10){\line(-1,0){0.2}}
\put(1.27,0.34){$\times$}
\put(2.27,3.31){$\times$}
\put(3.27,0.49){$\times$}
\put(4.27,1.75){$\times$}
\put(5.27,0.31){$\times$}
\put(6.27,0.69){$\times$}
\put(8.27,0.24){$\times$}
\put(9.27,1.33){$\times$}
\put(10.27,6.02){$\times$}
\put(11.27,0.64){$\times$}
\put(12.27,0.06){$\times$}
\put(13.27,8.39){$\times$}
\put(14.27,0.40){$\times$}
\put(15.27,2.33){$\times$}
\put(16.27,3.98){$\times$}
\put(17.27,0.01){$\times$}
\put(18.27,2.73){$\times$}
\put(19.27,0.98){$\times$}
\put(1,0){\line(0,1){0.47}}
\put(2,0){\line(0,1){0.47}}
\put(1,0.47){\line(1,0){1}}
\put(1.5,0.47){\line(0,1){2.97}}
\put(1.5,3.44){\line(1,0){3}}
\put(3,0){\line(0,1){0.62}}
\put(4,0){\line(0,1){0.62}}
\put(3,0.62){\line(1,0){1}}
\put(3.5,0.62){\line(0,1){1.26}}
\put(5,0){\line(0,1){0.44}}
\put(6,0){\line(0,1){0.44}}
\put(5,0.44){\line(1,0){1}}
\put(7,0){\line(0,1){0.82}}
\put(5.5,0.44){\line(0,1){0.38}}
\put(5.5,0.82){\line(1,0){1.5}}
\put(6.5,0.82){\line(0,1){1.06}}
\put(3.5,1.88){\line(1,0){3}}
\put(4.5,1.88){\line(0,1){1.56}}
\put(2.5,3.44){\line(0,1){7.56}}
\put(8,0){\line(0,1){0.37}}
\put(9,0){\line(0,1){0.37}}
\put(8,0.37){\line(1,0){1}}
\put(8.5,0.37){\line(0,1){1.09}}
\put(8.5,1.46){\line(1,0){1.5}}
\put(10,0){\line(0,1){1.46}}
\put(9.5,1.46){\line(0,1){4.69}}
\put(11,0){\line(0,1){0.77}}
\put(11,0.77){\line(1,0){1.5}}
\put(12,0){\line(0,1){0.19}}
\put(12,0.19){\line(1,0){1}}
\put(13,0){\line(0,1){0.19}}
\put(12.5,0.19){\line(0,1){0.58}}
\put(11.5,0.77){\line(0,1){5.38}}
\put(9.5,6.15){\line(1,0){2}}
\put(14,0){\line(0,1){0.53}}
\put(15,0){\line(0,1){0.53}}
\put(14,0.53){\line(1,0){1}}
\put(14.5,0.53){\line(0,1){1.93}}
\put(14.5,2.46){\line(1,0){1.5}}
\put(16,0){\line(0,1){2.46}}
\put(15.5,2.46){\line(0,1){1.65}}
\put(15.5,4.11){\line(1,0){3}}
\put(17,0){\line(0,1){0.14}}
\put(18,0){\line(0,1){0.14}}
\put(17,0.14){\line(1,0){1}}
\put(19,0){\line(0,1){1.11}}
\put(20,0){\line(0,1){1.11}}
\put(19,1.11){\line(1,0){1}}
\put(17.5,0.14){\line(0,1){2.72}}
\put(19.5,1.11){\line(0,1){1.75}}
\put(17.5,2.86){\line(1,0){2}}
\put(18.5,2.86){\line(0,1){1.25}}
\put(16.5,4.11){\line(0,1){4.41}}
\put(10.5,6.15){\line(0,1){2.37}}
\put(10.5,8.52){\line(1,0){6}}
\put(13.5,8.52){\line(0,1){2.48}}
\put(0.65,-0.57){-9}
\put(1.65,-0.57){-8}
\put(2.65,-0.57){-7}
\put(3.65,-0.57){-6}
\put(4.65,-0.57){-5}
\put(5.65,-0.57){-4}
\put(6.65,-0.57){-3}
\put(7.65,-0.57){-2}
\put(8.65,-0.57){-1}
\put(9.8,-0.57){0}
\put(10.83,-0.57){1}
\put(11.83,-0.57){2}
\put(12.83,-0.57){3}
\put(13.83,-0.57){4}
\put(14.83,-0.57){5}
\put(15.83,-0.57){6}
\put(16.83,-0.57){7}
\put(17.83,-0.57){8}
\put(18.83,-0.57){9}
\put(19.65,-0.57){10}
\end{picture}

\vspace{0.4in}
{\bf Figure 3.}
{\small 
A realization of part of $\A_\infty$, approximating the local structure 
of $\A_n$
for large $n$.
The 2 visible ancestral lineages diverged at around time $16$.
}

\vspace{0.1in}
\noindent
Proposition \ref{prop:cond_ptpr} and the calculation above easily
imply the first assertion below; the second assertion
follows from the mixture representation (\ref{Amix}), where
in this setting the mixing makes no difference.
\begin{cor}
\label{cor:local}
Let $t_n \to \infty$ and let
$U_n$ be uniform on $\{1,2,\ldots,n\}$ independent of $\A_{t_n,n}$.
Write
$\{(U_n + i + \sfrac{1}{2}, s_{U_n +i}), \ i \in \mathbf{Z}\}$
for the point process associated with the lineage tree
$\A_{t_n,n}$, centered at lineage $U_n$, where 
$s_j = 0$ for $j$ outside $[1,n]$.
Then as $n \to \infty$ this point process converges in
distribution to the point process 
$\{(i+\sfrac{1}{2}, \eta_i), \ i \in \mathbf{Z}\}$
defining $\A_\infty$.
The same result holds for $\An$.
\end{cor}
Less formally, the structure of $\A_\infty$ around lineage $0$
provides an asymptotic approximation to the structure of
$\An$ around a random lineage.

\subsection{Some local calculations}
\label{sec:slc}
We record some elementary calculations within $\A_\infty$,
reflecting aspects of the $n \to \infty$ behavior of the
lineage trees $\An$. For a lineage at time $s$ we call the present
(time - $0$) number of species descending from this lineage the 
{\em size} of this lineage.
The $n \to \infty$ limit of
$n^{-1}(\mbox{number of lineages in $\An$ at time $s$})$
is what we will call the {\em density of lineages} 
in $\A_\infty$ at time $s$.

\begin{cor}
\label{calc}
[Some calculations for $\A_\infty$.]
\\
 {\bf (a)} The density of ancestral lineages at time $s$ in the past
 equals $(1+s)^{-1}$, and the size 
of a random lineage
at time $s$
 has Geometric$((1+s)^{-1})$ distribution;\\
 {\bf (b)} the rate of lineages merging as $s$ increases (time runs
 backwards) is $m(s)=2(1+s)^{-1}$ , and given that this
 event occurs at $s$ for some lineage then the size of the lineage it merges with
has
 Geometric$((1+s)^{-1})$ distribution;\\
{\bf (c)} as $s$ decreases (time runs forward) the rate at which a lineage of size $k\geq 1$
 branches is $b_k(s)=(k-1)(s(1+s))^{-1}$ at time $s$, and the size of
 the lineage produced on the left of the branchpoint has Uniform
 distribution on  $\{1,\dots,k-1\}$.\\
\end{cor}
  
\begin{proof}
(a) The density of ancestral lineages
at time $s$ in the past 
is just the density of branching points at times greater
than $s$  
\[ G(s) = \int_s^\infty \ f(u) du = (1+s)^{-1}.\] 
Hence, the number of extant species descended from a ``typical''
lineage 
at time $s$ has Geometric($(1+s)^{-1}$) distribution
\[ p_s(i)= \left( \sfrac{1}{1+s} \right)\left(\sfrac{s}{1+s}
\right)^{i-1}, \ i \geq 1\]  
since this is the distribution of distances between branchpoints
at heights greater than $s$.

(b) As $s$ increases (time runs backwards) the
probability of a lineage merging with another lineage is 
\[ m(s)= 2 \sfrac{f(s)}{G(s)} = \sfrac{2}{1+s}\]
because such a merger occurs in $[s,s+ds]$ when
one of the two branchpoints separating the given lineage
from its neighboring lineages, which must be at height $ \geq s$,
occurs during $[s,s+ds]$, and this has chance 
$f(s) ds/G(s)$ for each branchpoint.
Moreover, if a lineage merges at $s$ then
(independent of the size of the first lineage)
the size of the second lineage has
Geometric${((1+s)^{-1})}$ distribution above.

(c) As $s$ decreases (time runs forwards), the unconditional rate
of mergers of clades of sizes $k_1,k_2$ at time $t$
(per unit time, relative to number of species) equals
\[ G(s) (1-G(s))^{k_1-1} f(s) (1-G(s))^{k_2-1} G(s) \]
which we observe by considering the required heights of branchpoints
for this event to occur.
Similarly the number of size $k_1+k_2$ lineages at time $t$,
relative to number of species, equals
\[ G(s) (1-G(s))^{k_1+k_2-1} G(s) . \]
Thus the rate of splitting of a size $k_1+k_2$ lineage into
two lineages of sizes $k_1, k_2$ equals
\[\frac{G(s) (1-G(s))^{k_1-1} f(s) (1-G(s))^{k_2-1} G(s) }
{ G(s) (1-G(s))^{k_1+k_2-1} G(s) }
= \frac{1}{s(1+s)}. \]
Thus, if a lineage is of size $k$ then at time $s$ the stochastic rate of
branching is
\[ b_k(s) = \sfrac{k-1}{s(1+s)}\]  
Since  the rate of splitting is independent of the choice of partition
of $k$ into $k_1$ and $k_2$ the size of a left
subclade lineage is Uniform on $\{1,2,\ldots,k-1\}$. 
\end{proof}

\section{Time reversal and consequences}
Recall that for a {\em stationary} Markov process, its
time-reversal is also a stationary Markov process.
For a Markov process which is not stationary, or which is conditioned
on a terminal value,
the time-reversal is typically {\em non}-homogeneous.
So Lemma \ref{lem:timerev} below highlights a special feature of our processes.

In the critical branching process underlying our model
(Section \ref{sec:model}),
the population size
is the continuous-time Markov chain with transition rates
\begin{equation}
q_{i,i+1} = q_{i,i-1} = i .
\label{CBP:rates}
\end{equation}
Recall the definition of the complete tree $\Tn$.
Write
$(N_n(s), \Tor_n \geq s \geq 0)$
for the associated process which counts the number of species
at time $s$ before present. 
The next lemma makes precise a sense in which the process
$(N_n(s))$ is a time-reversal of the chain (\ref{CBP:rates}) 
started at $0$.
\begin{lem}
\label{lem:timerev}
Let 
$(\Nhat_n(s), 0 \leq s \leq T_n^0)$
be the continuous-time chain (\ref{CBP:rates})
with $\Nhat_n(0) = n$, run until the first hitting
time $T_n^0$ on state $0$.  Then
\[ (N_n(s), \Tor_n \geq s \geq 0)
\ed
(\Nhat_n(s), 0 \leq s \leq T_n^0)
. \]
\end{lem}
\begin{proof}
We verify that $(\Nhat_n(s), 0 \leq s \leq T_n^0)$ is the time-reverse of
the population size process by checking probabilities of primitive
events (see Section \ref{sec:ext_anc} for more sophisticated views).
Fix
$s_M > s_{M-1} > \ldots > s_1 > s_0 = 0$
and positive integers
$1=k_M, k_{M-1},\ldots,k_1 = n$
with $|k_m - k_{m-1}| = 1$.
Set $k_{M+1} = 0$.
The event
\begin{quote}
as $s$ decreases, $N_n(s)$ jumps from 
$k_{m+1}$ to $k_m$ during $[s_m, s_m + ds_m]$
($\forall M \geq m \geq 1$)
and makes no other jumps
\end{quote}
has measure
\[ ds_M \times 
\prod_{m=M}^2 
\left( e^{-k_m(s_m - s_{m-1})} \ k_m \ ds_{m-1} \right) 
\times e^{-k_1s_1}
 \]
where the first term $ds_M$ comes from the uniform Bayes prior.
For the reversed process,
the event
\begin{quote}
as $s$ increases, $\Nhat_n(s)$ jumps from 
$k_{m}$ to $k_{m+1}$ during $[s_m, s_m + ds_m]$
($\forall 1 \leq m \leq M$)
and makes no other jumps
\end{quote}
has probability
\[
\prod_{m=1}^M 
\left( e^{-k_m(s_m - s_{m-1})} \ k_m \ ds_m \right) 
.  \]
By inspection, the first measure is exactly $1/n$ times the second probability,
so after conditioning the probability measures are equal.
\end{proof} 


We now observe two simple consequences of this time-reversal
identity.
The process
$(\Nhat_n(s), 0 \leq s \leq T_n^0)$
is a skip-free martingale started at $n$ and run until hitting
$0$, so by the hitting time formula for martingales
\[ \prob\left(\max_{0 \leq s \leq T_n^0} \Nhat_n(s) \geq c\right) = \sfrac{n}{c}, 
\quad c \geq n . \]
So Lemma \ref{lem:timerev} implies
\begin{cor}
\label{cor:max_popn}
\[ \prob\left(\max_{\Tor_n \geq s \geq 0} N_n(s) \geq c\right) = \sfrac{n}{c}, 
\quad c \geq n . \]
\end{cor}
Second, every extinction within the process $\Tn$ corresponds
to a downwards step in $N_n(s)$ as $s$ decreases and hence to an upwards
step in $\Nhat_n(s)$ as $s$ increases.
The number of such upward steps equals 
$(D_n - n)/2$,
where $D_n$ is the number of steps of the embedded jump chain
of $\Nhat_n(\cdot)$, which is just discrete-time simple symmetric
random walk.
\begin{cor}
\label{cor:conv_popn}
Within the model $\Tn$ of a clade on $n$ extant species,
the total number $N_n^{\mathrm{ext}}$ of extinct species is distributed as
$(D_n - n)/2$, where $D_n$ is the hitting time to $0$ for
simple symmetric random walk started at $n$.
In particular
\[ n^{-2} D_n \cd \sfrac{1}{2}\tau_1 \]
where $\tau_1$ is the first passage time of standard
Brownian motion from $1$ to $0$, with density function
  \[ f_{\tau_1}(x)=(2 \pi x^3)^{-\frac{1}{2}}\, e^{-\frac{1}{2x}},\quad
  0<x<\infty . \] 
\end{cor}
The second assertion follows, of course, from weak convergence
of simple random walk to Brownian motion.

\section{Exploiting the contour process}
\label{sec:contour}
The results so far answer some, but not all, questions
one might ask about the complete tree $\Tn$ and the lineage
tree $\An$ in our model.
For instance, the time-reversed process
$(\Nhat_n(s))$ in Lemma \ref{lem:timerev} 
has a $n \to \infty$ rescaled limit,
the well-known {\em Feller branching diffusion}, 
which therefore is the limit of the population size process
$(N_n(s), \Tor_n \geq s \geq 0)$.
But this doesn't tell us anything about
the relationship between 
$(N_n(s))$ and the lineage tree $\An$.
For instance, a conceptually interesting question in the species
context concerns 
$N_n(\Tmrca_n)$, 
the total number of species in the clade alive at the time of 
most recent common ancestor of the extant species.
Recall also that the results of Section \ref{sec:loc_anc}
were all based on the exact formula in Proposition \ref{prop:cond_ptpr},
but we have not yet given any indication of its proof.
It turns out that both these matters, and the local limit structure of the complete tree, can be studied using the
{\em contour process}, described next.

\subsection{The contour process}
For any deterministic population process in continuous time,
starting at the birth of one single individual,
in which individuals
have birth times, death times, and may give birth to children at distinct times,
there is a particular representation as a 
{\em rooted planar tree}
which we now describe.
Each individual is represented by an edge whose
length equals that individual's lifetime.
The birth of an offspring corresponds to a branchpoint from its parent's edge,
and the length of the parent's edge up to this
branchpoint equals the age of the parent at this offspring's birth time.
From the branchpoint the offspring's edge is drawn to the right of the
parent's edge.
If the total population is finite, then we can label the individuals
in a ``depth-first'' search order. This is
illustrated in Figure 4 where tree edges have been drawn as full vertical
lines and the branchpoints have been indicated by horizontal dotted lines.

\begin{figure}[hhhh]
\centerline{\epsfig{figure=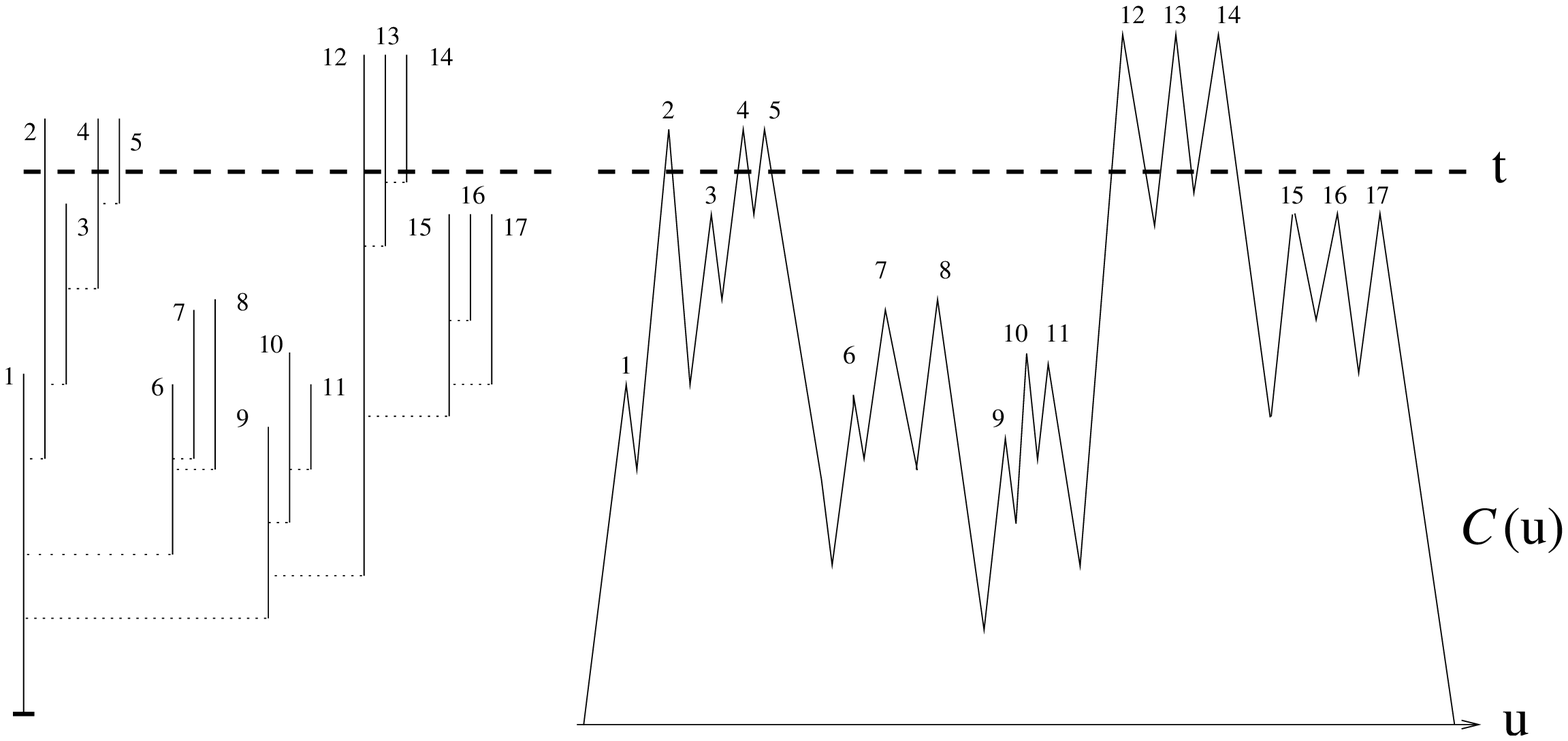,height=2.1in,width=4.25in,angle=0}}
\vspace{0.2in}
{\bf Figure 4.} A realization of a tree $\T_{t,n}$ with $n=6$ extant
individuals (labeled $\{2,4,5,12,13,14\}$) and its
contour process representation $\C(u)$.
\end{figure}

\noindent
Associated with such a rooted planar tree is its
{\em contour process} defined as follows
(these ideas go back to Neveu and Pitman
\cite{NP89a} 
and their broader significance can be seen in the lecture
notes of Pitman \cite{csp}).
The contour process $\C(u)$
is a continuous function giving the distance from the root at time
$u$ in a unit-speed depth-first walk around the tree.
Such a walk starts at the root, traverses each edge completely once
upwards and once downwards, following the depth-first order
(intuitively: clockwise around the edges of the tree)
ending back at the root.
So the contour process consists of alternating line segments
of slopes $+1$ (``rises") and slopes $-1$ (``falls").
The unit speed convention implies that heights in the contour process
match the times in the population process (birth and death times are
matched respectively by the local minima and local maxima in the
contour process). 

\subsection{Contour processes associated with random trees}
Recall that $\T$ denotes the continuous-time critical branching 
process started at time $0$ with one individual and continued
until extinction.
The relevance of contour processes is indicated in the next result of
Neveu-Pitman-Le Gall 
\cite{LG89,NP89a}.
\begin{prop}
\label{lem:contour} In the contour
process of $\T$ the sequence 
of rises and falls $(\xi_1,-\xi_2,\xi_3,-\xi_4,\ldots,\xi_{M-1})$,
excluding the last fall, has the distribution 
derived from a 
sequence $(\xi_i)_{i\geq 1}$ of independent Exponential$(1)$ variables, 
for 
$M:= \min \{m: \xi_1 - \xi_2 + \xi_3 - \xi_4 + \ldots - \xi_m < 0\}.$
\end{prop}
Call this contour process $(\xi_1, -\xi_2,\ldots,\xi_{M-1})$ an
{\em ERW excursion}, for Exponential random walk.
Accordingly call the infinite sequence 
$(\xi_1,-\xi_2,\xi_3, -\xi_4, \ldots)$ an {\em ERW process}.  
Here is a classical result

\begin{lem}
\label{lem:max_ht}
Let $H$ be the maximum height in an ERW excursion, or
equivalently (by Proposition \ref{lem:contour}) the extinction time
of $\T$.  Then
\[ \prob(H > h) = (1+h)^{-1}, \quad 0 < h < \infty . \]
\end{lem}

\begin{proof}
This follows directly from the law of the population size
process of $\T$ given in (\ref{Feller}). The extinction time of $\T$ is
greater than $h$ if and only if the population size of $\T$
at time $h$ is strictly greater than $0$. By (\ref{Feller}) the probability
of this is $1-h(1+h)^{-1}=(1+h)^{-1}$.  
\end{proof}

Before proceeding to new results let us indicate the
proof \cite{lea-paper} of Proposition \ref{prop:cond_ptpr},
because our arguments in subsequent sections will use similar ideas.
Fix $t$ and $n$.
Condition the contour process $\C(\cdot)$
to have exactly $n$ upcrossings over height $t$;  
see Figure 4.
This gives the contour process of the random tree
($\Ttn^+$, say)
which is the CBP conditioned on having exactly $n$ individuals
alive at time $t$.
This $\Ttn^+$ is the same as our model $\Ttn$ except for
the ``direction of time parameter" convention, and except for the
fact that in $\Ttn$ the process terminates with the $n$ individuals
at the present time, whereas in $\Ttn^+$ the process of descendants of
these $n$ individuals continues until extinction.
But the latter difference plays no role in the following argument.
The heights of the minima between each pair of
successive upcrossings in Figure 4 match the divergence of lineages
of that pair of extant individuals. Marking these heights
at regular horizontal interval spacings 
gives exactly the point process $\A_{t,n}$ as in Figure 1 
except for  reflecting the vertical time scale.
Since $\C(\cdot)$ is strong Markov and stationary the parts of an ERW excursion between a downcrossing of $t$
and the next upcrossing of $t$ are mutually independent, and moreover are distributed exactly as the  
reflection of the original ERW excursion that is conditioned not to have height greater than $t$.
Thus these heights of lineage divergence, when measured on the
reflected time scale (i.e. downwards from $t$), are
distributed as the maximum height $H$ in Lemma \ref{lem:max_ht}
conditioned on $\{H<t\}$.  
This conditioned distribution is the distribution
(\ref{def:fts}), as required for Proposition \ref{prop:cond_ptpr}.

\subsection{Species numbers at time of most recent common ancestor and
weak convergence of the contour process}

Recall that $N_n(\Tmrca_n)$ stands for the number of species alive at
the time of the most recent common ancestor. 
In the contour process the number of species at any time $s$
after its origin is the number of up-crossings (which equals the number of down-crossings) in the contour process
at height $s$. If the time since the origin of $\Tn$ is $\Tor_n=t$
then the contour process  has $n$ up- and down-crossings at
height $t$. If the time of the most recent common ancestor is
$\Tmrca_n=s$ then the maximal depth of the subexcursions below height
$t$, measured away from $t$, is $s$; see Figure 5. 
In other words, the lineage divergence of the most recent
common ancestor is the lowest local minimum between the first and last upcrossing of $t$ and occurs at height $t-s$. 

\begin{figure}[hhhh]
\centerline{\epsfig{figure=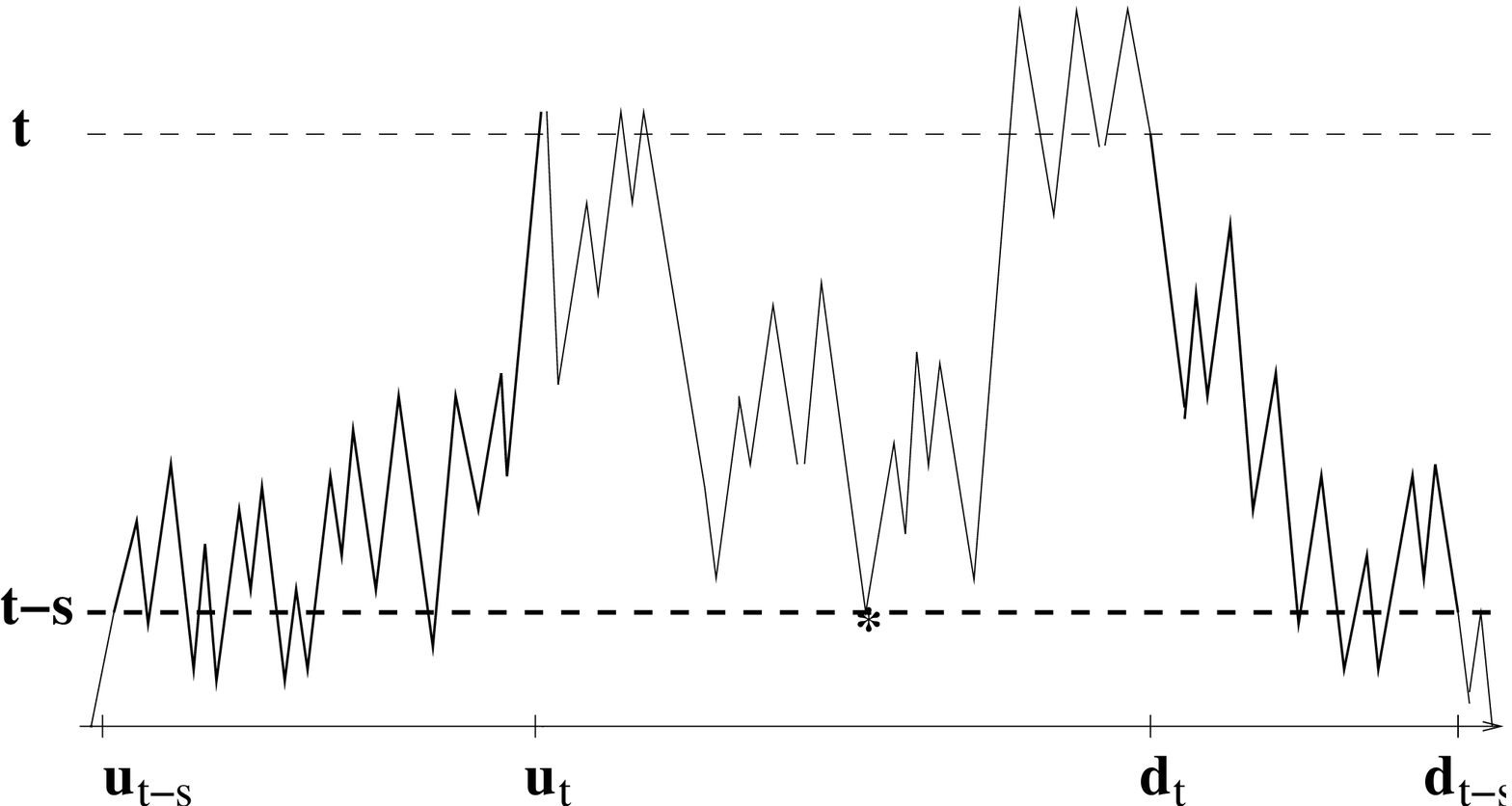,height=2.1in,width=4.25in,angle=0}}
\vspace{0.2in}
{\bf Figure 5.} {Parts of the contour process between $[u_{t-s},u_t]$
  and $[d_t,d_{t-s}]$ describe the number
of species alive at the time of the most recent common ancestor.}
\end{figure}

\noindent
In the contour process mark by  $u_s$ the horizontal coordinate of the
first upcrossing of 
a height $s$ and $d_s$ the coordinate of the last downcrossing of this height.
There are no up- or down-crossings of $t-s$ before the first
upcrossing of $t-s$ and after the last downcrossing of $t-s$.
And if $t-s$ is the height of the $\Tmrca$, as in Figure 5, then there are no up- or
downcrossings of $t-s$ between the first
upcrossing of $t$ and the last downcrossing of $t$. So, $N_n(\Tmrca_n)$
is the number of upcrossings of $t-s$ between $u_{t-s}$ and $u_t$, plus
the number of downcrossings between $d_t$ and $d_{t-s}$.
Since the contour process is an ERW excursion that is
conditioned to have $n$ upcrossings and downcrossings at height
$\Tor_n$, we can now calculate  

\begin{lem} 
\label{lem:num_mrca}
Conditional on $(\Tor_n,\Tmrca_n)=(t,s)$, $N_n(\Tmrca_n)$
is distributed as a sum of two independent Geometric($p_n$) random variables, where
\[p_n=1-\frac{t-s}{1+t-s}\,\frac{s}{1+s}.\] 
\end{lem}

\begin{proof}
Since the contour process $\C(\cdot)$ is strong Markov and stationary,
the part of the process between the first upcrossings $u_{t-s}$ of
$t-s$ and $u_t$ of $t$ when considered from height $t-s$ upwards:
$\;\C(u)-(t-s),\; u_{t-s}\leq u\leq u_t,\;$ is distributed as an ERW process  
conditioned to reach height $s$ before it reaches a depth $-(t-s)$
and stopped when it first hits $s$. Since $\C(\cdot)$ has the same
law when its $u$ coordinate is run in reverse,
the part of the contour process between the last downcrossings
$d_t$ of $t-s$ and $d_{t-s}$ of $t$ when run backwards in the $u$ coordinate:
$\;\C(u)-(t-s),\; d_{t-s}\geq u\geq d_t,\;$ is also distributed as a
ERW process  
conditioned to reach height $s$ before it reaches a depth $-(t-s)$ and
stopped when it first hits $s$. Additionally these two parts of the
contour process are independent. 

The probability an ERW process reaches $s$ before it reaches $-(t-s)$,
by the  law of maximum height $H$ in  Lemma \ref{lem:max_ht}, is
\[
\frac{\prob(H>s)}{\prob(H>s)+\prob(H>t-s)-\prob(H>s)\prob(H>t-s)}=\frac{1+t-s}{1+t}.
\]
The probability an ERW process makes $k$ upcrossings of $0$ until it first
hits $s$ , provided its
height stays below $s$  and its depth above $-(t-s)$, is for $k=1,2,\dots$
\[
\left(\prob(H<t-s)\prob(H<s)\right)^{k-1}\prob(H>s)=
\left(\frac{t-s}{1+t-s}\frac{s}{1+s}\right)^{k-1}\frac{1}{1+s},
\]
So the number of upcrossings of $t-s$ $\C(u)$ makes during
$[u_{t-s},u_t]$  
has a Geometric$\left(1-\sfrac{t-s}{1+t-s}\sfrac{s}{1+s}\right)$
distribution.
\end{proof}

Since by Corollary \ref{cor:t_or,mrca} $(n^{-1}\Tor_n,n^{-1}\Tmrca_n)\cd (\Tor,\Tmrca)$, as $n \to \infty$
\[
np_n =
1-\frac{\Tor_n-\Tmrca_n}{1+\Tor_n-\Tmrca_n}\frac{\Tmrca_n}{1+\Tmrca_n}
\cd
\frac{1}{\Tor-\Tmrca}+\frac{1}{\Tmrca}
\]
and the above two Geometric($p_n$) variables, when rescaled by
$n^{-1}$, converge to independent 
Exponential$(\lambda(\Tor,\Tmrca))$ variables, where
$\lambda(t,s)=(t-s)^{-1}+s^{-1}$. Consequently,
the conditional law of $n^{-1}N_n(\Tmrca_n)$ given $(\Tor_n,\Tmrca_n)$
converges
to a Gamma variable with shape parameter $2$ and scale
parameter $\lambda(\Tor,\Tmrca)$. Combining this with the result of
Corollary \ref{cor:t_or,mrca} we have 
established assertion (\ref{ftriple}) below.
 
\begin{cor} 
\label{cor:conv_triple}
The joint limit behavior of the triple $\Tor_n,\Tmrca_n,N_n(\Tmrca)$
is given by 
\[
(n^{-1}\Tor_n,n^{-1}\Tmrca,n^{-1}N_n(\Tmrca_n)) \cd
(\Tor,\Tmrca,\Nmrca)
\]
where the limit has the joint density 
\begin{eqnarray}
f_{\Tor,\Tmrca,\Nmrca}(t,s,r)&=&t^{-2}s^{-2}\lambda(t,s)^2re^{-\frac{1}{s}-\lambda(t,s)
  r}\nonumber\\
&=&(t-s)^{-2}s^{-4}re^{-\frac{1}{s}-\frac{tr}{s(t-s)}}
,\quad 0<s<t,0<r.
\label{ftriple}
\end{eqnarray}
The marginal density of $\Nmrca$ is
\[
f_{\Nmrca}(r)=2(1+r)^{-3},\quad r>0.
\] 
\end{cor}

\noindent 
The marginal density formula follows from (\ref{ftriple})
via a calculus exercise.
Note that while the distribution of $\Nmrca$ has mean $1$ it has infinite variance.

{\bf Remark.}
The contour process of $\Ttnn$
(illustrated in Figure 5),
in the limit $t_n/n \to t\in(0,\infty)$,
converges after rescaling to a Brownian excursion, conditioned
on total local time at height $t$ being equal to $1$.
Results like Corollary \ref{cor:conv_triple} may be reinterpreted as
providing exact formulas for quantities defined in terms of
such conditioned Brownian excursions.

\subsection{Extinct species}\label{sec:ext_anc}
Textbooks (e.g. \cite{PH-ME} page 24) often say
\begin{quote}
the probability that a given fossil is actually part of an
ancestral lineage [of some extant species] is actually rather remote.
\end{quote}
Various calculations relevant to this issue can be done within our model.

Consider some species $v$ that originated at time $h$ before the present.
In the limit as $n\to\infty$, the distribution of the clade of species
descending from $v$ is given by 
the local limit structure of the complete tree $\T_\infty$. 
As stated in Section \ref{sec:local_complete}, in the 
limit,  the descendants of a 
species $v$ evolve, as time runs forwards, as in an ordinary critical
branching process $\T$. 
Then, the chance that some descendant of $v$ (or $v$ itself) is extant
at present equals the 
chance of the survival of its descendant tree $\T$ for time $h$
or longer. By Lemma \ref{lem:max_ht} this is precisely $(1+h)^{-1}$, so we have

\begin{cor}
\label{cor:ext_t}
For any species alive at time $h$ before the present, the chance that
some of its descendant species (or the species itself) is extant is,
in the limit $n\to\infty$, 
\[
1/(1+h) .
\]
\end{cor} 

Now consider the total number $\Nanc_n$ of species
that are ancestral to the $n$ extant ones. 
(Precisely, we exclude the extant species, and go back to the
time of origin of the clade).
Intuitively, because (Lemma \ref{lem:timerev})
the number of species at time $h$ is $N_n(h) \approx n$ for $h = o(n)$,
and because (Corollary \ref{cor:t_or,mrca}) the time of origin
$\Tor_n$ is of order $O(n)$, we expect from Corollary
\ref{cor:ext_t} that
\[ \Ex [
\Nanc_n ] \approx 
\int_0^{O(n)} \sfrac{n}{1+h} \ dh
\approx n \log n . \]
We shall prove a precise result as Corollary \ref{cor:conv_anc},
based on the following lemma.

\begin{lem}
\label{lem:anc}
Conditional on $\Tor_n = t$,
the total number of ancestral individuals $\Nanc_n$ in $\Tn$
satisfies 
\[\Nanc_n \mathop{=}^d \mathop{\sum}\limits_{i=1}^{n} X_i\]
where $X_i, 1\leq i\leq n$ are independent, $X_1$ has  
Poisson($t$) distribution
and $X_2,\dots,
X_n$ have the law, with $f_t(\cdot)$ as in (\ref{def:fts}),
\[
\prob(X_i=k)=\int_0^t\frac{e^{-s}s^k}{k!}f_t(s)ds,\quad
k\geq 0.
\]
\end{lem}  

\begin{proof}
Label the extant individuals $\{1,2,\dots,n\}$ from left to right 
as they appear in the contour process. Let $X_i$ be the
number of ancestors of the $i$th extant individual, without
including any of the ones previously counted in $X_j, j<i$.

Suppose $\Tor=t$, then the ancestry of the extant individuals is
described by the part of the contour process $\C(\cdot)$ below height $t$.
Recall that the part of $\C(\cdot)$ below $t$ consists
of: $n-1$ independent sub-excursions below $t$, which we label $e_i, 1\leq i\leq
n-1$, and
the part of $\C(\cdot)$ before the first up-crossing and after the last
down-crossing of $t$; we label the former part as $e_{0,R}$.  See
Figure 6.

Let $h_i, 1\leq i\leq n-1$ be the depths of the sub-excursions $e_i$,
so that $t-h_i$ are the heights of the lowest points of $e_i$. 
These match the times of divergence of lineages of  extant
individuals. Their law was given by (\ref{def:fts}) of Proposition
\ref{prop:cond_ptpr}. Now, partition the excursions $e_i$ at their
lowest points and let $e_{i,R}, 1\leq i\leq n-1$ denote the parts on the right.
Then Figure 6 shows that the ancestors of the $1$st extant
individual correspond in the piece of
the contour process $e_{0,R}$ to the  levels of constancy of the process 
$\varsigma_{0,R}(u)=\mathop{\inf}_{v\geq u}(e_{0,R}(v))$. These levels
of constancy of  $\varsigma_{0,R}$ match the times of lineage
divergence of the ancestors of individual $1$.
Similarly for the $i$th, $2\leq i\leq n$, extant individual Figure 6 shows
that the additional ancestors of individual $i$ (excluding those appearing as ancestors of extant
individuals $j<i$) correspond in
the piece of the contour process $e_{i-1,R}$ to the levels of
constancy of the process $\varsigma_{i-1,R}(u)=\mathop{\inf}_{v\geq
  u}(e_{i-1,R}(v))$.  

\begin{figure}[hhhh]
\centerline{\epsfig{figure=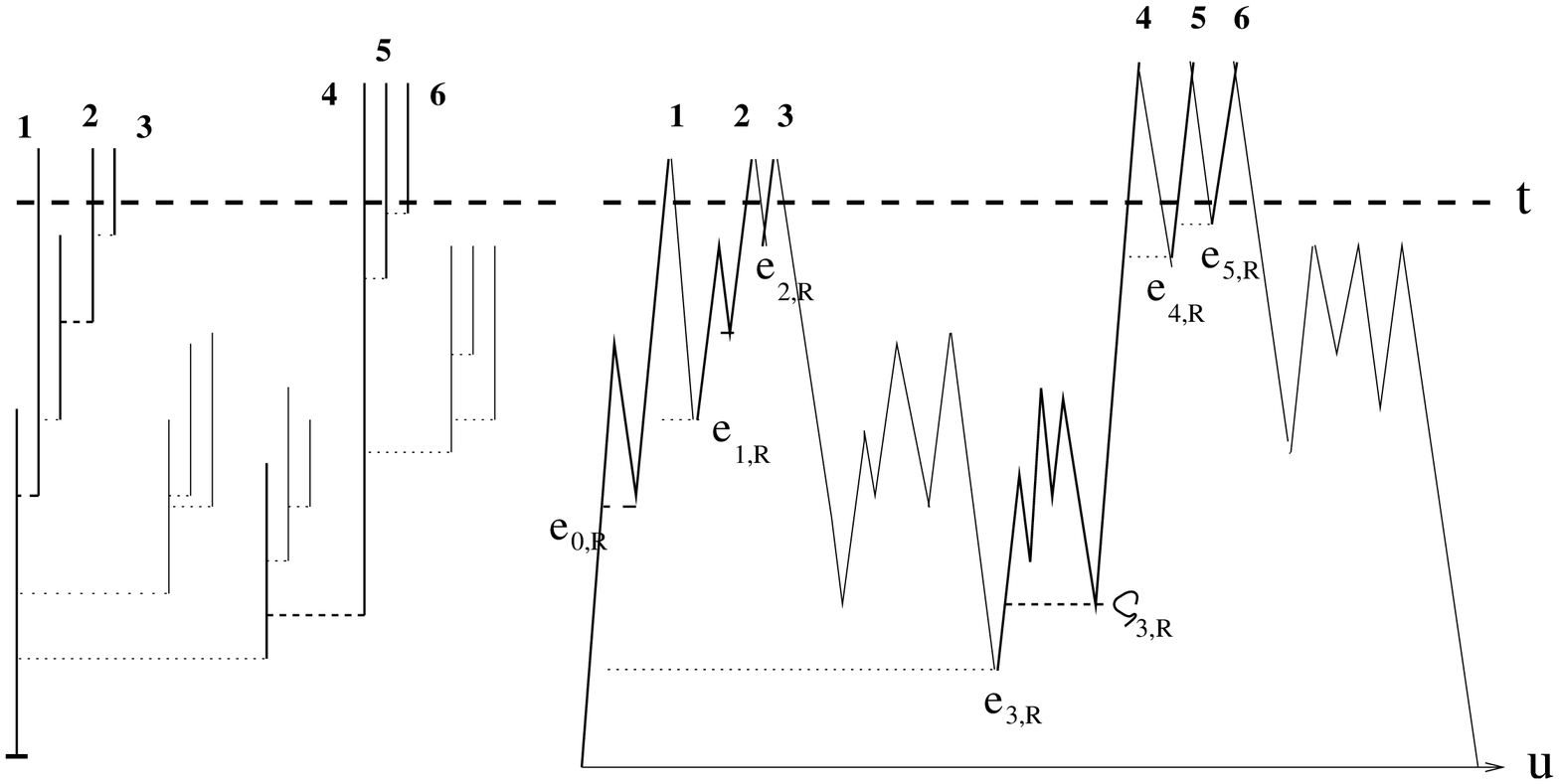,height=2.1in,width=4.25in,angle=0}}
\vspace{0.2in}
{\bf Figure 6.} {Ancestral lineages of the extant individuals
  (labeled $\{1,2,3,4,5,6\}$) are matched in the contour process by the levels
  of constancy of the processes $\varsigma_{i-1,R}$, for $1\leq i\leq n$.}
\end{figure}

\noindent
So the number of ancestors $X_i$ of the $i$th extant individual is the
number of levels of constancy 
of the process $\varsigma_{i-1,R}(\cdot)$.
It is clear that the piece $e_{0,R}$ is distributed as an ERW process
conditioned to hit $t$ before $0$, and stopped the first time it hits $t$. 
It is less obvious but none the less true (see Lemma~6 of
\cite{lea-paper}) that, given $h_i$, the piece $e_{i,R}$ is also
distributed as an ERW process conditioned to reach $h_i$ before $0$,
and stopped the first time it hits $h_i$. For such a conditioned ERW, 
the levels of constancy of its future infimum process 
form a Poisson process, restricted to, respectively, the set $[0,t]$
for $e_{0,R}$,  and
$[0,h_i]$ for $e_{i-1,R},\, 2\leq i\leq n$. 
This can easily be seen for levels of constancy of the past supremum
process of a conditioned ERW, (Lemma~6 of \cite{lea-paper}), then time
reversibility of an ERW excursion implies the rest.
So the number of ancestors of the $1$st extant individual is
Poisson$(t)$, and the number of additional ancestors of the extant
individuals $i$, $2\leq i\leq n$, is Poisson$(h_i)$. Combining this with
the distributions of the depths $h_i$, given  by (\ref{def:fts}) in Proposition
\ref{prop:cond_ptpr}, 
we have proved the claim.
\end{proof}

For the limit of the number of ancestors $\Nanc_n$ we have
\begin{cor}
\label{cor:conv_anc} 
As $n \to \infty \ $ $\frac{\Nanc_n}{n\log{n}} 
\cp 1$.
\end{cor}
\begin{proof}
Fix $(t_n)$ such that $t_n/n \to t \in (0,\infty)$.
Because (Corollary \ref{cor:t_or,mrca}) 
$n^{-1} \Tor_n$ has a distributional limit on $(0,\infty)$,
it suffices to prove the following:
conditional on $\{\Tor_n = t_n\}$
we have
$\frac{\Nanc_n}{n \log n} \cp 1$.

We shall prove this using the representation
$\Nanc_n=\sum_{i=1}^n X_i$ 
from Lemma \ref{lem:anc}, where in the following argument 
we are always 
conditioning on $\{\Tor_n = t_n\}$.
Note that the contribution to the sum from $X_1$ is
negligible
(because $X_1$ has Poisson($t_n$) distribution),
so we may assume $X_1$ has the same distribution as the 
$X_i$, $2 \leq i \leq n$.
We now calculate
\[ \Ex [X_2] = 
\int_0^{t_n} s f_{t_n}(s) \ ds
\sim \int_0^{t_n} s(1+s)^{-2} \ ds
\sim \log t_n
\sim \log n \]
and a similar calculation shows
\[ \var [X_2] = O(n) . \]
Thus
\[ E [\Nanc_n] \sim n \log n ; \quad 
\var [\Nanc_n] = O(n^2) \]
and the desired result
$\frac{\Nanc_n}{n \log n} \cp 1$
follows via Chebyshev's inequality.
\end{proof}

\subsection{Local limit structure of the complete
  tree}\label{sec:local_complete} 
The contour process makes it conceptually easy to see a result,
complementing Corollary \ref{cor:local} for the lineage tree $\A_n$,
concerning the local limit behavior of the complete tree 
$\Tn$. It turns out that the local structure
relative to a given typical individual in $\Tn$, converges to the local
structure relative to the root in an infinite tree that can be easily
defined from a CBP tree. 
There are two versions of such results, depending on
whether for the typical individual we choose 
a random {\em extant} species or a random species from the entire
history of the clade (which will be {\em extinct}, with
probability $\to 1$ as $n \to \infty$). The details in obtaining the
results are rather fussy, so we only outline the proofs. 

Let $i$ be an individual in the complete tree $\Tn$, 
with birth time $b(i)$ say. Within this section our convention for the
time parameter in $\Tn$ is that it increases as time increases. 
For $\sigma > 0$ 
let $\widetilde{\T}_n(i,[b(i) - \sigma, b(i)+\sigma])$
denote the subtree of $\Tn$ comprised of all the individuals $j$ whose
birth time is in the time interval $[b(i) -\sigma, b(i)+\sigma]$
and for whom the divergence time of their lineage from that of $i$
is in the time interval $[b(i) - \sigma, b(i)+\sigma]$.
See Figure 7.
Call $i$ the 
{\em distinguished individual} in
$\widetilde{\T}_n(i,[b(i) -\sigma,b(i)+\sigma])$.

\setlength{\unitlength}{0.25in}
\begin{picture}(15,15)(0,-2)
\put(0,-1){\vector(0,1){12}}
\put(0,0){\line(-1,0){0.2}}
\put(0.1,-0.1){$b(i) - \sigma$}
\put(0,5){\line(-1,0){0.2}}
\put(0.1,4.9){$b(i) $}
\put(0,10){\line(-1,0){0.2}}
\put(0.1,9.9){$b(i) + \sigma $}
\put(-1.3,10.5){time}
\put(5,1){\line(0,1){5.3}}
\put(6,4){\line(0,1){1.7}}
\multiput(5.2,3.9)(0.5,0){2}{-}
\put(7,5){\line(0,1){3}}
\multiput(6.2,4.9)(0.5,0){2}{-}
\put(6.6,6.5){$i$}
\put(8,7.5){\line(0,1){2}}
\multiput(7.2,7.4)(0.5,0){2}{-}
\put(9,6.2){\line(0,1){2.7}}
\multiput(8.2,6.1)(0.5,0){2}{-}
\multiput(7.2,6.1)(0.5,0){2}{-}
\put(10,8.4){\line(0,1){2.9}}
\multiput(9.2,8.3)(0.5,0){2}{-}
\put(11,7.1){\line(0,1){3.6}}
\multiput(10.2,7.0)(0.5,0){2}{-}
\multiput(9.2,7.0)(0.5,0){2}{-}
\put(12,2.4){\line(0,1){3.1}}
\multiput(5.2,2.3)(0.5,0){14}{-}
\put(13,4.5){\line(0,1){3.8}}
\multiput(12.2,4.4)(0.5,0){2}{-}
\put(14,7.7){\line(0,1){2.9}}
\multiput(13.2,7.6)(0.5,0){2}{-}
\put(15,5.5){\line(0,1){0.8}}
\multiput(13.2,5.4)(0.5,0){4}{-}
\multiput(3,0)(0.6,0){22}{\_}
\multiput(3,10)(0.6,0){22}{\_}
\end{picture}

{\bf Figure 7.}
{\small 
The local structure of the complete tree, relative to
individual $i$.
}

\vspace{0.12in}
\noindent
We now describe
an infinite random tree $\widetilde{\T}$ derived from the CBP.
Take a distinguished individual, born at time $0$, 
and let the tree of it and its descendants be distributed as
the CBP tree $\T$.
Let the parent of this individual have 
Exponential($1$) age at time $0$
and have an independent Exponential($1$) lifetime after 
time $0$.
Inductively, let the grandparent have
Exponential($1$) age at 
the birth of the parent, 
and have an independent Exponential($1$) lifetime after 
that birth time; and so on.
For each of these ancestors, let them have other children during their
lifetimes at the times of a Poisson (rate $1$) process and
let the trees of such children and their descendants be 
distributed as independent CBP trees. 
This completes the description of $\widetilde{\T}$.

Recall the construction (Proposition \ref{lem:contour})
of CBP tree $\T$ from the ERW excursion 
$(\xi_1,-\xi_2,\xi_3,\ldots)$.
It is easy to check that given a two-sided 
ERW process
$(\ldots, -\xi_{-2}, \xi_{-1},-\xi_0,\xi_1,-\xi_2,\xi_3, \ldots)$
an analogous construction produces the infinite tree
$\widetilde{\T}$.
Write
$\widetilde{\T}[-\sigma, \sigma]$
for the subtree of $\T$ comprised of individuals $j$ whose
birth time is in the time interval $[-\sigma, \sigma]$
and for whom the divergence time of the lineages of $j$  
and of the distinguished individual
is in the time interval $[ - \sigma, \sigma]$.
Note that 
\begin{equation}
\mbox{
$\widetilde{\T}[-\sigma,\sigma]$
 is determined by 
$(\xi_i, M^- \leq i \leq M^+)$}
\label{M+-}
\end{equation}
where 
$(\xi_{M^-}, - \xi_{M^- +1},\ldots, -\xi_0, \xi_1, \ldots , \xi_{M^+ -1}, - \xi_{M^+})$
is the excursion of the two-sided ERW process above height $- \sigma$.

Here is the result for the convergence of the local structure of $\Tn$
as seen relative to a random (extinct) individual, to that of
the local structure of $\widetilde{\T}$.
\begin{prop}\label{thm4}
Let $I_n$ denote a uniform random species from the clade $\Tn$.
Then as $n \to \infty$ for fixed $\sigma > 0$,
\[ \widetilde{\T}_n(I_n,[b(I_n)-\sigma, b(i) + \sigma])
\cd \widetilde{\T}[-\sigma,\sigma] . \]
\end{prop}
{\bf Remark.}
The underlying notion of {\em convergence} of finite trees
is the natural one, which can be formalized in several equivalent ways,
e.g. via a point process representation.

\begin{proof}
We outline the proof, omitting details.
Write $(\xi_i, i \geq 1)$
for the ERW process.
Fix an integer $m \geq 2$.
Let $\theta_{m,N}$ be the empirical distribution
of the $N$ $2m$-vectors
$\{(\xi_{2i+1},\xi_{2i+2},\ldots,\xi_{2i+2m}); \ 0 \leq i \leq N-1\}$.
So $\theta_{m,N}$ is a random probability distribution, 
which (using the Glivenko-Cantelli Theorem on $\mathrm{R}^{2m}$)
converges in probability, as $N\to\infty$, to the non-random probability
distribution $\mu_m = \mathrm{dist}(\xi_1,\ldots,\xi_{2m})$.
By large deviation theory
(see \cite{DZ92} 
\S 6.3)
this convergence remains true conditional on events
$A_N$ for which $1/P(A_N) = O(\beta^N)$ for all $\beta > 1$.

To prove the proposition, recall (Lemma \ref{lem:t_or})
that $\Tor_n$ is order $n$.
So we can fix $(t_n)$ such that $t_n/n \to t \in (0,\infty)$
and it is sufficient to prove the Proposition for
$\Ttnn$.
Fix also integers
$N_n$ such that $N_n/n^2 \to v \in (0,\infty)$.
Let $A_{N_n}$ be the event that an ERW process has an excursion above
$0$ with  exactly $N_n$ rises and falls, and that this excursion
has exactly $n$ upcrossings over height $t_n$.
(Then $n^2$ is precisely the right scaling for 
the number of rises and falls of an excursion with $n$ upcrossings of
a level $t_n$ of  order $n$.) 
Conditioned on this event, the ERW excursion
is the contour process of a random tree
$\TtnnN^+$ which is the tree $\Ttnn$ continued until extinction
that is conditioned to have the total number of individuals equal to $N_n$ .
Let us first prove the proposition for $\TtnnN^+$.

One can show that 
the probability $P(A_{N_n})$ decreases not faster than polynomially
in $1/{N_n}$, so by our ``large deviation" result earlier, 
the empirical distribution $\theta_{m,N_n}$ of $2m$-tuples
conditioned on $A_{N_n}$ converges to $\mu_m$.  
This implies the weaker result that, for 
$J_n$ uniform on 
$\{2,4,6,\ldots,2N_n\}$, 
\begin{equation}
 (\xi_{J_n -m+1},\ldots,\xi_{J_n},\ldots,\xi_{J_n +m})
\cd \mu_m \label{Jn}
\end{equation}
where the left side is conditioned on $A_{N_n}$.
But this says that, relative to a uniform random individual
$I_n$ in $\TtnnN^+$,
any aspect of 
the ``local structure" of the tree which is determined
by the contour process segment of length $2m$ centered 
on that individual will converge in distribution to
the same aspect of the local structure of 
$\widetilde{\T}$.
By taking $m$ large and appealing to (\ref{M+-}),
we see that the proposition holds for 
$\TtnnN^+$.

To complete the proof it is enough to show that the proposition
holds for the stopped tree $\TtnnN$.
Unfortunately this does not follow directly from the unstopped case,
because a non-negligible fraction of all individuals in
$\TtnnN^+$ will be descendants of the $n$ individuals alive at
time $t_n$ after origin.
Instead, fix small $0 < \delta_1 < \delta_2$
and consider the segments of the contour process $\C^+$ of
$\TtnnN^+$ defined by:\\
\indent $s_1$ is the segment of $\C^+$ until its first upcrossing of 
$(1-\delta_1)t_n$,\\
\indent $s_2$ is the segment of $\C^+$ from the subsequent downcrossing of
$(1-\delta_2)t_n$ until the next upcrossing of 
$(1-\delta_1)t_n$,\\
\indent $s_3$ is the segment of $\C^+$ from the subsequent downcrossing of
$(1-\delta_2)t_n$ until the next upcrossing of 
$(1-\delta_1)t_n$,\\
\indent \ldots\\
\indent $s_N$ is the final segment of $\C^+$ after the final downcrossing of $t_n$.

\noindent Conditional on the event $A_{N_n}$, 
there is some conditional distribution of starting and ending
positions for each segment.
Given all these positions, 
each segment is distributed as an ERW process conditioned on having the first
upcrossing of a certain level after a prescribed number of steps.   
The number of these segments is stochastically bounded as $n \to
\infty$, so the probability of the conditioning event for each segment
is still only polynomially small in $1/$length of the segment. 
Thus separately on each segment we can show as above that the
contour process satisfies (\ref{Jn}) for $J_n$ uniform on that
segment.
Since these segments comprise (in the $n \to \infty$ limit) a proportion
$1 - \eps(\delta_1,\delta_2)$ of the entire contour process
of $\TtnnN$,
where $\eps \to 0$ as $\delta_1, \delta_2 \to 0$,
we can deduce the proposition 
for the stopped process $\TtnnN$.
\end{proof}

We now state (omitting the similar argument) the parallel local limit result for $\Tn$
as seen from a random {\em extant} individual.
In this setting the relevant limit infinite tree,
which we again call $\widetilde{\T}$,
is a variation of the $\widetilde{\T}$ above described as follows.
The distinguished individual has 
Exponential($1$) age at time $0$.
Its ancestors and their descendants are all as described before, except
that now the infinite tree $\widetilde{\T}$ is stopped at time $0$.
\begin{prop}\label{thm5}
Let $I_n$ denote a uniform random extant species from the clade $\Tn$.
Then as $n \to \infty$ for fixed $\sigma > 0$,
\[ \widetilde{\T}_n(I_n,[-\sigma, 0])
\cd \widetilde{\T}[-\sigma,0] . \]
\end{prop}

One can now make exact calculations of probabilities for the distinguished
individual in $\widetilde{\T}$,
which represent the $n \to \infty$ limit results for a random
extant individual in $\Tn$.
Here is a simple example of possible calculations within $\widetilde{\T}$.
\begin{cor}\label{conseq6} 
For the distinguished individual in $\widetilde{\T}$:\\
  {\bf (a)} the probability that its parent is extant equals $1/2$;\\
  {\bf (b)} the probability that some ancestor of it is extant equals
  $1-e^{-1}$.
\end{cor}
\begin{proof}
(a) The probability that the parent of the distinguished individual is alive at time $0$ is simply
$\prob(\xi_1 < \xi_2)$, where $\xi_1$ is the age
of the distinguished individual, and $\xi_2$ is the subsequent lifetime
of its parent after the birth.  Because $\xi_1$ and $\xi_2$ are independent exponential$(1)$
random times, we have $\prob(\xi_1 < \xi_2)=1/2$, by symmetry. 

(b) To calculate the probability that no ancestor of the distinguished individual is
still alive, one only need to note that the times at which some
ancestor originates form a Poisson process 
of rate $1$, and an ancestor originating at time $s$ before present has chance
$e^{-s}$ to be extant, so the random number of extant ancestors
has Poisson distribution with mean
$\int_0^\infty e^{-s} \times 1 \ ds = 1$,
and thus takes value $0$ with probability $e^{-1}$.
\end{proof}

\section{Final remarks }\label{sec:final}

{\bf 1.}
Our model of $\Tn$ and $\An$ has considerable variability between
realizations.
This can be seen mathematically in our distributional formulas
(Corollary \ref{cor:conv_triple} in particular) 
and visually on
our web site \cite{me-Phylo-site}. 
In one sense this variability is an artifact of the uniform prior on
time of origin, but serves a useful purpose in emphasizing that
radically different appearance of real-world trees might logically be just
chance variation without biological significance.

{\bf 2.}
Wollenberg et al 
\cite{WAA96} study via simulation a model similar to ours -- critical
branching conditioned on $n$ extant species -- but handle the issue
of time of origin in a different way, by taking it as the
deterministic time $t_n$ which is the maximum likelihood estimator
of origin time.
In a sense this is unrealistic in the opposite sense to that of
the previous remark, by underestimating variability.
Our model extends more naturally to higher-order taxa.

{\bf 3.}
Our model is qualitatively similar 
(in the sense of orders of magnitude)
to the Moran model, for quantities which can be studied in the
latter model.
In fact the results involving local weak limits
(sections \ref{sec:llp} and \ref{sec:slc})
are exactly the same in our model
as in a continuized Moran model,
because our model converges 
(in the $n \to \infty$ limit)
to the continuized Moran model over time intervals
(backwards from present)
of length $t = o(n)$.

{\bf 4.}
Neutral models like ours are unrealistic for large clades,
by the following reasoning.
For an $n$-species clade, our model gives
(Corollary \ref{cor:t_or,mrca})
the time of origin of a clade as order $n$ time units ago.
The time unit is mean species lifetime, 
typically estimated as a few million years.
Thus our model predicts the origin of a $n$-species clade
as being at least $n$ million years ago, which is known to
be an overestimate for most clades of size $n \geq 100$.

{\bf 5.}
The local point process limit in Corollary \ref{cor:local}
is a simple instance of a general notion of
{\em local weak convergence} of graphical structures
associated with point processes on $\mathbf{R}^d$ or abstract spaces.
See \cite{me101,me110} for more sophisticated examples.
In particular, 
Proposition \ref{thm4} fits the general setting of
{\em asymptotic fringe distributions} which exist for many
different models of random trees \cite{me52}.

{\bf 6.}
Mathematicians traditionally tend to regard pictures as mere visual aids to illustrate a logical argument.  
But the {\em graphical representations} we use
in Figures 1 and 4 really comprise the essence of the mathematical
argument, by relating our model of random trees to well-understood
models of point processes or random walks.

\paragraph{Acknowledgement}
We thank Maxim Krikun for helpful comments.


\def\cprime{$'$} \def\polhk#1{\setbox0=\hbox{#1}{\ooalign{\hidewidth
  \lower1.5ex\hbox{`}\hidewidth\crcr\unhbox0}}} \def\cprime{$'$}
  \def\cprime{$'$} \def\cprime{$'$}
  \def\polhk#1{\setbox0=\hbox{#1}{\ooalign{\hidewidth
  \lower1.5ex\hbox{`}\hidewidth\crcr\unhbox0}}} \def\cprime{$'$}
  \def\cprime{$'$} \def\polhk#1{\setbox0=\hbox{#1}{\ooalign{\hidewidth
  \lower1.5ex\hbox{`}\hidewidth\crcr\unhbox0}}} \def\cprime{$'$}
  \def\cprime{$'$} \def\cydot{\leavevmode\raise.4ex\hbox{.}} \def\cprime{$'$}
  \def\cprime{$'$} \def\cprime{$'$} \def\cprime{$'$} \def\cprime{$'$}
  \def\cprime{$'$} \def\cprime{$'$} \def\cprime{$'$} \def\cprime{$'$}
  \def\cprime{$'$}

\end{document}